
\def\LaTeX{%
  \let\Begin\begin
  \let\End\end
  \let\salta\relax
  \let\finqui\relax
  \let\futuro\relax}

\def\UK{\def\our{our}\let\sz s}
\def\USA{\def\our{or}\let\sz z}

\UK



\LaTeX

\USA


\salta

\documentclass[twoside,11pt]{article}
\setlength{\textheight}{24cm}
\setlength{\textwidth}{17cm}
\setlength{\oddsidemargin}{2mm}
\setlength{\evensidemargin}{2mm}
\setlength{\topmargin}{-15mm}
\parskip1mm


\usepackage[usenames,dvipsnames]{color}
\usepackage{amsmath}
\usepackage{amsthm}
\usepackage{amssymb}
\usepackage{bbm}
\usepackage{bm}
\usepackage[mathcal]{euscript}
\usepackage{hyperref}
\usepackage{enumitem}
\allowdisplaybreaks[4]

%
%


\definecolor{viola}{rgb}{0.3,0,0.7}
\definecolor{ciclamino}{rgb}{0.5,0,0.5}
\definecolor{rosso}{rgb}{0.85,0,0}

\def\comm #1{{\color{blue}#1}}

\def\comm #1{#1}




\bibliographystyle{plain}


%
\newtheorem{theorem}{Theorem}[section]
\newtheorem{remark}[theorem]{Remark}
\newtheorem{corollary}[theorem]{Corollary}

\newtheorem{proposition}[theorem]{Proposition}
\newtheorem{lemma}[theorem]{Lemma}

\finqui

\def\Bcenter{\Begin{center}}
\def\Ecenter{\End{center}}
\let\non\nonumber




\def\step #1 \par{\medskip\noindent{\bf #1.}\quad}


\def\Lip{Lipschitz}
\def\Holder{H\"older}
\def\Fre{Fr\'echet}
\def\Lady{Lady\v zhenskaya}

\def\rhs{right-hand side}



\def\multibold #1{\def\arg{#1}%
  \ifx\arg\pto \let\next\relax
  \else
  \def\next{\expandafter
    \def\csname #1#1#1\endcsname{{\boldsymbol #1}}%
    \multibold}%
  \fi \next}

\def\pto{.}

\def\multical #1{\def\arg{#1}%
  \ifx\arg\pto \let\next\relax
  \else
  \def\next{\expandafter
    \def\csname #1#1\endcsname{{\cal #1}}%
    \multical}%
  \fi \next}


\def\multimathop #1 {\def\arg{#1}%
  \ifx\arg\pto \let\next\relax
  \else
  \def\next{\expandafter
    \def\csname #1\endcsname{\mathop{\rm #1}\nolimits}%
    \multimathop}%
  \fi \next}

\multibold
qwertyuiopasdfghjklzxcvbnmQWERTYUIOPASDFGHJKLZXCVBNM.

\multical
QWERTYUIOPASDFGHJKLZXCVBNM.


\multimathop
diag dist div dom mean meas sign supp .


\def\Accorpa #1#2 #3 {\gdef #1{\eqref{#2}--\eqref{#3}}%
  \wlog{}\wlog{\string #1 -> #2 - #3}\wlog{}}
\def\Accorparef #1#2 #3 {\gdef #1{\ref{#2}--\ref{#3}}%
  \wlog{}\wlog{\string #1 -> #2 - #3}\wlog{}}


\def\<#1>{\mathopen\langle #1\mathclose\rangle}
\def\norma #1{\mathopen \| #1\mathclose \|}

\def\[#1]{\mathopen\langle\!\langle #1\mathclose\rangle\!\rangle}

\def\ioT {\int_0^T}

\def\intQ{\int_Q}
\def\iO{\int_\Omega}

\def\dt{\partial_t}
\def\dn{\partial_\nnn}

\def\checkmmode #1{\relax\ifmmode\hbox{#1}\else{#1}\fi}

\def\aeQ{\checkmmode{a.e.\ in~$Q$}}
\def\aeG{\checkmmode{a.e.\ on~$\Gamma$}}


\def\erre{{\mathbb{R}}}
\def\enne{{\mathbb{N}}}




\def\genspazio #1#2#3#4#5{#1^{#2}(#5,#4;#3)}
\def\spazio #1#2#3{\genspazio {#1}{#2}{#3}T0}

\def\L {\spazio L}
\def\H {\spazio H}
\def\W {\spazio W}

\def\C #1#2{C^{#1}([0,T];#2)}


\def\Lx #1{L^{#1}(\Omega)}
\def\Hx #1{H^{#1}(\Omega)}

\def\Hdn{H^2_{\nnn}(\Omega)}



\let\theta\vartheta

\let\eps\varepsilon
\let\ph\varphi

\let\TeXchi\chi                         
\newbox\chibox
\setbox0 \hbox{\mathsurround0pt $\TeXchi$}
\setbox\chibox \hbox{\raise\dp0 \box 0 }
\def\chi{\copy\chibox}



\let\emb\hookrightarrow

\def\cd{C_\delta}

\def\muh{\hat\mu}

\def\umin{\bu_{\min}}
\def\umax{\bu_{\max}}
\def\Uad{{{\cal U}_{\mathrm{ad}}}}
\def\UR{{{\cal U}_{R}}}

\def\Vp{V^*}
\def\Vsp{\VVV^{*}_\sigma}

\def\ov #1{{\overline{#1}}}

\def\bv{{\boldsymbol{v}}}
\def\bu{{\boldsymbol{u}}}
\def\opt{\widetilde{\boldsymbol{u}}}
\def\bh{{\boldsymbol{h}}}
\def\bph{{\widetilde{\ph}}}
\def\bw{{\boldsymbol{w}}}
\def\bs{{\boldsymbol{\sigma}}}
\def\0{{\boldsymbol{0}}}

\def\Sol{{\cal S}}
\def\Jred{{\cal J}_{\rm red}}
\def\UR{{\cal U}_R}
\def\YY{{\cal Y}}

\def\ZZ{{\cal Z}}
\def\rem{{\cal R}}

\def\bvh{{\bv^\bh}}
\def\phh{{\ph^\bh}}
\def\muh{{\mu^\bh}}
\def\omh{{\omega^\bh}}

\def\bz{{\boldsymbol{\zeta}}}
\def\bva{{\bv^a}}
\def\pa{{p^a}}
\def\pha{{\ph^a}}
\def\mua{{\mu^a}}
\def\oma{{\omega^a}}

\def\Hs{{\HHH_\sigma}}
\def\Vs{{\VVV_\sigma}}
\def\Ws{{\WWW_\sigma}}

\def\Hqn{{H^4_{\nnn}(\Omega)}}

\def\alphaphQ{{\alpha_1}}
\def\alphaphO{{\alpha_2}}
\def\alphacon{{\alpha_3}}

\def\solst{(\bv,p,\ph,\mu,\omega)}
\def\soln{(\bv_n,p_n,\ph_n,\mu_n,\omega_n)}
\def\solopt{(\widetilde{\bv},\widetilde{p},\widetilde{\ph},\widetilde{\mu},\widetilde{\omega})}
\def\solsth{(\bvh,p^\bh,\phh,\muh,\omh)}
\def\sollin{(\bw,q,\psi,\theta,\w)}
\def\soladj{(\bva,\pa,\pha,\mua,\oma)}

\def\diff{(\phh-\ph)}

\def\w{\mathrm{w}}

\Begin{document}


%

\title{
Optimal Control of a Navier--Stokes--Cahn--Hilliard System \\
for Membrane-fluid Interaction
}

\date{}
\author{}
\maketitle
\Bcenter
\vskip-1cm
{\large\sc Andrea Signori$^{(1)}$}\\
{\normalsize Email: {\tt andrea.signori@polimi.it}}\\[.25cm]
{\large\sc Hao Wu$^{(2)}$}\\
{\normalsize Email: {\tt haowufd@fudan.edu.cn}}\\[.5cm]
$^{(1)}$
{\small Dipartimento di Matematica, Politecnico di Milano}\\
{\small via E. Bonardi 9, 20133 Milano, Italy}
\\
{\small Alexander von Humboldt Research Fellow}
\\[.3cm]
$^{({2})}$ 
{\small School of Mathematical Sciences, Fudan University}
\\
{\small Handan Road 220, 200433 Shanghai, China}
%
%
%
%
%

\Ecenter

\Begin{abstract}
\noindent
We consider an optimal control problem for a two-dimensional Navier--Stokes--Cahn--Hilliard system arising in the modeling of fluid-membrane interaction. The fluid dynamics is governed by the incompressible Navier--Stokes equations, which are nonlinearly coupled with a sixth-order Cahn--Hilliard type equation representing the deformation of a flexible membrane through a phase-field variable.
Building on the previously established existence and uniqueness of global strong solutions for the coupled system, we introduce an external forcing term acting on the fluid as the control variable. Then we seek to minimize a tracking-type cost functional, demonstrating the existence of an optimal control and deriving the associated first-order necessary optimality conditions. A key issue is to establish sufficient regularity for solutions of the adjoint system, which is crucial for the rigorous derivation of optimality conditions in the fluid dynamic setting.


\End{abstract}

\salta
\pagestyle{myheadings}
\newcommand\testopari{\sc Signori \ -- \ Wu}
\newcommand\testodispari{{\sc Optimal control of a Navier--Stokes--Cahn--Hilliard system}}
\markboth{\testodispari}{\testopari}
\finqui


\section{Introduction}
\label{SEC:INTRO}

The investigation of membrane formation and its dynamic behavior represents a central topic of interdisciplinary research that intersects the fields of physics, biology, and engineering. Biological membranes, for instance, perform critical roles by delineating cellular interiors from the extracellular environment and partitioning intracellular spaces into distinct functional compartments. Their geometry, which often becomes highly complex in response to mechanical and physicochemical stimuli, is closely tied to their biological performance~\cite{SeifertLipowsky1995}.  Similarly, in engineering applications such as polymer electrolyte membrane fuel cells, a hydrophobic polymer membrane separates the anode and cathode, facilitating the selective transport of ion and electron~\cite{PromislowWetton2009}.
In recent decades, numerous studies have been aimed at deciphering the complex morphological transformations and multiscale features of membranes. These efforts include modeling and numerical simulations~\cite{
Aland2017, AlandEgererLowengrubVoigt2014,
CampeloHernandezMachado2006,
Canham1970,DuLiuRyhamWang2005, DuLiuRyhamWang2005b, DuLiuRyhamWang2009, DuLiuWang2004, DuLiuWang2006, Helfrich1973, LowengrubRatzVoigt2009}
as well as rigorous mathematical analyses~\cite{ChengWangWiseYuan2020, ClimentEzquerraGuillenGonzalez2019,  ColliLaurencot2012, DuLiLiu2007,EntringerBoldrini2015, LiuTakahashiTucsnak2012, WuXu2013}.

At the molecular level, membranes mainly consist of amphiphilic bilayers. The classical Helfrich model describes the equilibrium shape of membranes via the bending elasticity energy, see, e.g., \cite{Canham1970, Helfrich1973}.
In its simplest isotropic form, it can be expressed as
\begin{align*}
	{{\cal E}}_{\text{elastic}} = \frac{k}{2}  \int_{\Gamma} (H - H_0)^2 \, \mathrm{d}S,
\end{align*}
where $\Gamma$ denotes the evolving membrane,
$k$ denotes the bending rigidity, $H$ is the mean curvature of the membrane and $H_0$ represents the spontaneous curvature. For homogeneous bilayer lipid membranes, spontaneous curvature $H_0$ usually vanishes.
The complex interplay between stresses arising from membrane elasticity and (external) hydrodynamic flows can lead to a rich phenomenology. To capture the non-equilibrium evolution and deformation of membranes influenced by the surrounding fluid, it is crucial to incorporate hydrodynamic coupling effects.
In this regard, the phase-field (or diffuse-interface) method offers a straightforward and efficient framework to capture and integrate a wide range of physical phenomena on multiple scales~\cite{DuLiuRyhamWang2005, LowengrubRatzVoigt2009, YueFengLiuShen2004}.
In the phase-field framework, a phase-field variable $\ph$ is introduced to distinguish the interior and exterior regions of the membrane by the values $+1$ and $-1$, respectively. The Helfrich bending energy, with $H_0=0$ for simplicity, can be approximated by a modified Willmore energy functional:
\begin{align}\label{Willmore}
	{{\cal E}}_\varepsilon(\ph) = \frac{k}{2\varepsilon} \int_{\Omega} \left(-\varepsilon \Delta \ph + \frac{1}{\varepsilon}(\ph^2 - 1) \ph \right)^2 \mathrm{d}x,
\end{align}
where $\varepsilon$ is a small positive parameter related to the interfacial thickness~\cite{DuLiuRyhamWang2005,DuLiuRyhamWang2009}. As $\varepsilon \to 0$, the phase-field model converges to the corresponding sharp-interface description~\cite{DuLiuRyhamWang2005}.
We postpone further technical details and references to the existing literature to a later stage, in order to first clearly present the main objectives of this contribution.

In this work, we address an optimal control problem for a hydrodynamic phase-field model that describes the migration and deformation of membranes in an incompressible viscous fluid. For simplicity, we consider the interior and exterior liquids of the membrane as Newtonian fluids with matched  densities (assumed to be one). Let $\Omega\subset \mathbb{R}^2$ be a bounded domain with smooth boundary $\Gamma:=\partial\Omega$ and let $T>0$ be the prescribed final time.
The control problem under consideration is governed by the following system of partial differential equations, which we will refer to as the \textit{state system}:
\begin{alignat}{2}
	\label{eq:1}
	& \dt \bv
	+ (\bv \cdot \nabla ) \bv
	- \div (2 \nu (\ph) D \bv)
	+ \nabla p
	=
	\mu \nabla \ph
	+ \bu
	\qquad && \text{in $Q$},
	\\
	\label{eq:2}
	& \div \bv =0
	\qquad && \text{in $Q$},
	\\
	\label{eq:3}
	& \dt \ph + \bv \cdot \nabla \ph
	- \div (m(\ph)\nabla \mu)
	=0
	\qquad && \text{in $Q$},
	\\
	\label{eq:4}
	& \mu = - \eps\Delta \omega
	+   \dfrac{1}{\eps} f'(\ph) \omega
	+ \eta \omega
	\qquad && \text{in $Q$},
	\\
	\label{eq:5}
	& \omega =
	-\eps \Delta \ph
	+  \dfrac{1}{\eps} f(\ph)
	\qquad && \text{in $Q$},
		\\
	\label{eq:6}
	& \bv = \0,
	\quad
	\dn \ph = \dn \Delta \ph = m(\ph)\dn \mu =0
	\qquad && \text{on $\Sigma$},
	\\
	\label{eq:7}
	& \bv(0)=\bv_0,
	\quad
	\ph(0)=\ph_0
	\qquad && \text{in $\Omega$}.
\end{alignat}
\Accorpa\Sys {eq:1} {eq:7}
Here, we set $Q:= \Omega\times (0,T)$ and $\Sigma:= \Gamma\times (0,T)$.
In \eqref{eq:6}, the vector $\nnn$ denotes the outward unit normal field on the boundary $\Gamma$ and $\dn$ stands for the corresponding outward normal derivative.
The system \Sys\ can be understood as a coupling between the Navier--Stokes equations for the velocity field $\bv: \Omega\times (0,T)\to \mathbb{R}^2$ and a six-order convective variant of the classical fourth-order Cahn--Hilliard equation for the phase-field function $\varphi: \Omega\times (0,T)\to \mathbb{R}$.
The vector $\bu$ appearing on the right-hand side of \eqref{eq:1} acts as a distributed forcing term and will later serve as a control variable. The tensor $D\bv$ denotes the symmetric part of the velocity gradient, defined as $D\bv = \frac{1}{2}(\nabla \bv + \nabla^{\top} \bv)$, while $p:\Omega\times(0,T)\to \mathbb{R}$ indicates the fluid's pressure. Moreover, the nonlinear term $\mu \nabla \ph$ models the Korteweg force associated with capillarity effects.

Alongside the state system \Sys, we now formulate the control problem, which consists in minimizing the tracking-type \textit{cost functional}:
\begin{align}
	\label{def:cost}
	{\cal J} (\ph,\bu)= \frac {\alphaphQ}2 \intQ|\ph-\ph_Q|^2\,\mathrm{d}x\mathrm{d}t
	+\frac {\alphaphO}2 \iO|\ph(T)-\ph_\Omega|^2\,\mathrm{d}x
	+\frac {\alphacon}2 \intQ|\bu|^2\,\mathrm{d}x\mathrm{d}t,
\end{align}
where the coefficients ${\alpha_i}$, $i=1,2,3$, are nonnegative real constants (with at least one $\alpha_i$ strictly positive to avoid a trivial scenario), and $\ph_Q:\Omega\times (0,T)\to \mathbb{R}$ and $\ph_\Omega:\Omega\to \mathbb{R}$ are prescribed target functions.
This minimization problem is subject to two constraints: the control variable $\bu$ should belong to a suitable admissible set $\Uad$ (specified below), and the state variables $(\bv,p,\varphi)$ must solve the associated state system \Sys\ corresponding to the given control $\bu$. The purpose of this study fits naturally within the analytical framework described above. Due to the need for analytical properties to ensure well-posedness (note that the system under investigation contains the Navier--Stokes equations as a subsystem), we focus on the two-dimensional case, where the existence and uniqueness of global strong solutions (see Theorem \ref{THM:RES} below) can be established following the argument in \cite{WY}.
For the control space $\mathcal{U}= \LLL^\infty(Q)= L^\infty(Q)^2$, we consider the \textit{admissible set} of control given by
\begin{align*}
	\Uad = \big\{
	\bu \in  \mathcal{U}  \,\,:\,\, u_{{\rm min}, i}\leq u_{ i}	\leq u_{{\rm max}, i} \quad \aeQ, \quad i=1,2
	\big\},
\end{align*}
for given and compatible fixed functions $\umin,\umax \in {\cal U}$.
Then, the optimal control problem that we are going to address reads as follows:
\begin{align*}
	\textbf{(CP)}
	\quad
	\min_{\bu \in \Uad} {\cal J} (\bu,\ph)
	\quad
	\text{subject to  $\solst$ that solves \Sys.}
\end{align*}
With regard to the problem \textbf{(CP)}, we shall prove the existence of optimal controls and characterize such optimal solutions by deriving the associated first-order necessary conditions. A detailed description of our results can be found in Section \ref{SEC:RES}.

We now provide some further explanations of the Navier--Stokes--Cahn--Hilliard system \Sys. The first two equations \eqref{eq:1}--\eqref{eq:2} describe the evolution of the divergence-free velocity field $\bv$ for the viscous incompressible fluid. Equations \eqref{eq:3}--\eqref{eq:5} represent a convective Cahn--Hilliard type equation, with the usual normalization for the phase-field variable $\ph$, denoting the difference of local concentrations. In the phase-field description, $\ph$ distinguishes the membrane interior ($\{\ph=1\}$) from the exterior ($\{\ph=-1\}$), and the membrane is modeled by a narrow, diffuse layer ($\{-1 < \ph < 1\}$) within a $\eps$-tubular neighborhood with $0<\eps\ll 1$. The function $m(\ph)$ in \eqref{eq:3} represents the diffusion mobility, which is assumed to be positive and nondegenerate. Define the total free energy functional (cf. \eqref{Willmore})
\begin{align}
 {\cal E}(\ph)
 & :=
 {\cal F}(\ph) + \eta {\cal G}(\ph)
 =
 \frac 12 \iO  \Big(-\eps\Delta\ph + \frac 1\eps f(\ph)\Big)^2\,\mathrm{d}x
 + \eta  \iO  \, \Big(\frac \eps2|\nabla\ph|^2 +\frac 1\eps F(\ph)\Big)\,\mathrm{d}x.
  \label{defE}
\end{align}
The variable $\mu$ in equation~\eqref{eq:4} is known as the {\it chemical potential} and corresponds to the first variation $\mu = \frac {\delta {\cal E}}{\delta \ph}$. Similarly, the variable $\omega$ denotes the first variational derivative of the Ginzburg--Landau free energy ${\cal G}$ such that $\omega = \frac {\delta {\cal G}}{\delta \ph}$.
In~\eqref{defE}, $F$ is a double-well shaped nonlinear potential, while $f$ denotes its derivative $f=F'$.
A canonical example of $F$ is the classical regular potential given by
\begin{align}
F(s)=\frac14(s^2-1)^2,\quad \forall\, s\in \erre.
\label{double-well}
\end{align}
It is important to mention that the total free energy ${\cal E}$ represents a higher-order extension of the Ginzburg--Landau free energy ${\cal G}$, which is intrinsically linked {with the classical} Cahn--Hilliard framework.
The parameter $\eta \in \erre$ controls the relative contribution of the Ginzburg--Landau free energy, that is, its sign and magnitude will influence the balance between the interfacial energy and the higher-order regularization effects.
When $\eta = 0$, the energy functional ${\cal E}$ reduces to the classical Willmore functional \eqref{Willmore} in the phase-field framework (see, e.g., \cite{ColliLaurencot2012,DuLiuWang2004}). In this case, the model accurately reproduces the Canham--Helfrich bending energy, which characterizes the elastic behavior of biomembranes based on their curvature. This connection has been extensively studied, for example, in \cite{DuLiuRyhamWang2005, DuLiuWang2004}, where the diffuse interface formulation is shown to converge to geometric models that involve the Willmore energy.
When $\eta > 0$, the energy functional ${\cal E}$ can be interpreted as a Willmore-type regularization of the Ginzburg--Landau free energy ${\cal G}$. The inclusion of a positive higher-order term enhances the regularity of the interface and penalizes the curvature concentration, making the model suitable for describing more complicated cases with anisotropic effects. Applications of this regularized energy include the modeling of thin film growth and coarsening in anisotropic systems, as discussed in \cite{BCMS,TorabiLowengrubVoigtWise2009}.
When $\eta < 0$, the energy functional ${\cal E}$ takes the form of the so-called functionalized Cahn--Hilliard (FCH) free energy, initially introduced in \cite{GompperSchick1990, Helfrich1973} to model phase separation in amphiphilic mixtures. Later, it was extended to capture nanoscale morphological changes in functionalized polymers and bilayer interfaces \cite{DaiPromislow2013, GavishHayrapetyanPromislowYang2011, PromislowWetton2009}. These soft materials, often membranes with charged groups, interact with polar solvents to form surfactant-rich interfaces.
Physically, the negative parameter $\eta$ reflects a competition between the square of the variational derivative of the classical Cahn--Hilliard free energy $\cal G$ and the energy itself \cite{DaiPromislow2013, GavishJonesXuChristliebPromislow2012}, balancing elastic deformation and hydrophilic surface effects. Unlike classical models that minimize the surface area, the FCH energy encourages an increase in the interfacial area, leading to complex structures such as stable bilayers and homoclinic interfaces \cite{DaiPromislow2013}.
This rich behavior has attracted extensive analytical and numerical studies, focusing on minimization, pattern formation, pearling instabilities, and topological transitions \cite{DaiPromislow2015, GompperSchick1990, PromislowWetton2009}. For further details, see \cite{DHPW, DaiPromislow2015, PQ, PromislowWu2017} and the references therein.

Returning to biological applications, the behavior of membranes plays a crucial role in many cellular processes. Given that membranes are typically closed surfaces, it is natural to incorporate physically relevant constraints such as volume conservation and surface area preservation. These constraints reflect the slow volumetric changes of the membrane and the inextensibility induced by amphiphilic molecules~\cite{AlandEgererLowengrubVoigt2014, DuLiuRyhamWang2009}. Within the diffuse-interface framework, these conditions are often approximated by suitable integral quantities:
\begin{align*}
A(\ph) = \int_\Omega \ph \, \mathrm{d}x,
\quad B(\ph) = \int_\Omega \left( \frac{\varepsilon}{2}|\nabla \ph|^2 + \frac{1}{\varepsilon}F(\ph)\right) \mathrm{d}x,
\end{align*}
where $F$ is given by \eqref{double-well}.
Enforcing the volume and surface area constraints can be realized either through Lagrange multipliers, or penalty terms added to the energy functional ${\cal F}$, see, e.g., \cite{DuLiuRyhamWang2005, DuLiuRyhamWang2005b, DuLiuWang2006}:
\begin{align}
{{\cal E}}_{\text{penalty}}(\ph) = {{\cal F}}(\ph) + \frac{M_1}{2}(A(\ph) - \beta)^2 + \frac{M_2}{2}(B(\ph) - \gamma)^2,
\label{energy-Pen}
\end{align}
where $\beta$ and $\gamma$ correspond to given values of volume and surface area, and $M_1, M_2 > 0$ are penalty parameters.
Although extensive analytical and numerical work has been devoted to phase-field membrane models without fluid coupling~\cite{CampeloHernandezMachado2006, ColliLaurencot2011, ColliLaurencot2012, DuLiuRyhamWang2005, DuLiuRyhamWang2005b, DuLiuWang2004, Wang2008}, the incorporation of fluid-membrane interaction significantly complicates the problem. Several hydrodynamic phase-field models, often derived via energetic variational principles, couple membrane elasticity with fluid dynamics \cite{AlandEgererLowengrubVoigt2014, DuLiuRyhamWang2009, GigaKirshteinLiu2018, LowengrubAllardAland2016}. In \cite{ClimentEzquerraGuillenGonzalez2016,ClimentEzquerraGuillenGonzalez2019}, the authors considered a Cahn--Hilliard--Navier--Stokes vesicle-fluid interaction model subject to the free energy ${{\cal E}}_{\text{penalty}}$ with $M_1=0$. In the Cahn--Hilliard framework, the volume
constraint is automatically guaranteed under the homogeneous Neumann boundary for the chemical potential (cf. \eqref{eq:3}, \eqref{eq:6}), and only the surface area constraint needs to be approximated by penalization.

When fluid coupling is neglected, that is, \(\bv = \mathbf{0}\) in \eqref{eq:3} and the momentum equation is absent, the analysis of sixth-order Cahn--Hilliard-type equations becomes more tractable and has attracted significant interest due to its wide range of applications \cite{CGSS6, CGSS8, SchimpernaWu2020, WY}. These include modeling oil-water-surfactant mixtures \cite{PZ1, PZ2, SP}, studying surface faceting phenomena \cite{KNR, KR}, and investigating phase-field crystal models \cite{GW1, GW2, M1, M2, WW}.
In the fluid-free setting, the FCH equation subject to periodic boundary conditions admits global weak solutions provided the potential is regular and the mobility may be degenerate, as rigorously established in \cite{DaiLiuLuongPromislow2021,DaiLiuPromislow2021}.
More recently, \cite{SchimpernaWu2020} extended the mathematical analysis to the physically significant logarithmic potential
\begin{align*}
	F_{\rm log}(s) = \frac{1}{2}(1+s)\ln(1+s) + \frac{1}{2}(1-s)\ln(1-s) - \frac{\theta}{2} s^2, \quad s \in (-1,1),
\end{align*}
overcoming technical challenges related to singular diffusion terms.
We also mention the recent work \cite{CGhyp}, in which the authors analyzed a corresponding hyperbolic relaxation (with regular potentials).
With regard to applications in optimal control problems, there are still very few works available. In this direction, we refer to \cite{CGSS6} that analyzes a control problem for a sixth-order Cahn--Hilliard system without fluid coupling, and \cite{ColliGilardiSprekels2018, FRS, PS, RoccaSprekels2015, SprekelsWu2021}, where velocity control problems are studied for hydrodynamic coupled systems involving the classical local or nonlocal Cahn--Hilliard equations.

We conclude by providing a brief overview of this work. In Section \ref{SEC:RES}, we introduce the notation and present the main results. Section \ref{SEC:DIFF} is dedicated to a thorough investigation of the differentiability properties associated with the solution operator of the system \Sys, laying the groundwork for our subsequent exploration of the first-order necessary conditions for optimality in the control problem. Finally, in Section \ref{SEC:CONTROL}, we address the control problem \textbf{(CP)}, where we not only demonstrate the existence of an optimal control, but also derive necessary optimality conditions utilizing appropriate adjoint variables. To derive these results, it was crucial to show that the hydrodynamic structure of the adjoint problem enables for sufficient regularity properties.

\section{Main Results}
\label{SEC:RES}
\setcounter{equation}{0}
\subsection{Notation}
\label{SUBSEC:NOT}
Let us first introduce some notation and conventions.
For any (real) Banach space $X$, we denote by $ \|\cdot\|_X $, $X^*$, and $\< \cdot,\cdot >_X$ the corresponding norm, its dual space, and the related duality pairing between $X^*$ and $X$. Given an interval $I\subset\mathbb{R}^+$, we introduce the function space $L^q(I;X)$ with $q\in [1,+\infty]$, which consists of Bochner measurable $q$-integrable functions (or essentially bounded if $q=+\infty$) with values in $X$. For $q\in [1,+\infty]$, $W^{1,q}(I;X)$ denotes the space of functions $f$ such that $f\in L^q(I;X)$ with $\partial_t f\in L^q(I;X)$, where $\partial_t $ means the vector-valued distributional derivative of $f$. For $q=2$, we set $H^1(I;X):=W^{1,q}(I;X)$.

Let $\Omega\subset \mathbb{R}^2$ be a bounded domain with a smooth boundary $\Gamma:=\partial\Omega$ and let $T>0$ be the prescribed final time. We set $Q:=\Omega \times (0,T)$, $\Sigma:=\Gamma \times (0,T)$ as well as
\begin{align*}
    Q_t & := \Omega \times (0,t),\quad \Sigma_t:=\Gamma \times (0,t), \quad \text{for $t\in(0,T)$}.
\end{align*}
The Lebesgue measure of $\Omega$ is denoted by $|\Omega|$. For any $q\in [1,+\infty]$, $k\in \mathbb{N}$, the standard Lebesgue and Sobolev spaces on $\Omega$ are denoted by $L^q(\Omega)$ and $ W^{k,q}(\Omega)$, respectively. When $k=0$, we have $W^{0,q}(\Omega)=L^{q}(\Omega)$ with the corresponding norm denoted by $\norma{\cdot}_{q}$. When $q = 2$, the convention $H^{k}(\Omega):=W^{k,2}(\Omega)$ will be used. As customary, we use bold letters to indicate vector- or matrix-valued spaces, e.g., $\LLL^q(\Omega) := L^q(\Omega; \erre^2)$, or  $\LLL^q(\Omega) := L^q(\Omega; \erre^{2\times2})$, for $q\in[1,+\infty]$. The same symbol will be used for the norm in some space and that in any power thereof.
For convenience, we introduce the notation
\begin{align*}
     &H:= \Lx2, \qquad V := \Hx1.
\end{align*}
The norm and inner product for $H$ are denoted by $\norma{\cdot}$ and $(\cdot,\cdot)_H$, respectively.
For every $f\in V^*$, we denote its generalized mean on $\Omega$ by
$\overline{f}=|\Omega|^{-1}\langle f,1\rangle_{V}$; if $f\in L^1(\Omega)$, then it holds $\overline{f}=|\Omega|^{-1}\int_\Omega f \,\mathrm{d}x$.
The linear subspace of $H$ with zero mean will be denoted by $H_{0}:=\big\{f\in H\ |\ \overline{f} =0\big\}$. We further define $V_0:=V\cap H_0$ and
\begin{align*}
     &\Hdn:= \{f \in H^2(\Omega): \dn f = 0 \,\,\aeG\, \},\\
     &\Hqn:= \{f \in H^4(\Omega): \dn f =  \dn \Delta f =0 \,\,\aeG\, \}.
\end{align*}
In connection with the spaces above, we adopt the usual framework of Hilbert triplets by identifying $H$ and $\Vp$ with subsets of $\Vp$ and $\Hdn^*$ in the usual way, that is,
\begin{alignat*}{2}
    \langle u,v \rangle_V & = \int_{\Omega} uv\,\mathrm{d}x,
     \qquad
    && \forall\,u\in H, \,\, v\in V,
    \\
    \langle u,w  \rangle_{\Hdn} & = \langle u ,w  \rangle_V,
    \qquad
    && \forall\,u \in \Vp, \,\,w \in \Hdn.
\end{alignat*}
Moreover, it holds $\Hdn \emb V \emb H \emb \Vp \emb \Hdn^*$ with dense and compact embeddings.
We also recall the following \Lady's inequality and Agmon's inequality
in two dimensions:
\begin{alignat}{2}
	\label{inter}
	\norma {v}_4 &\leq C \norma{v}^{\frac12}\norma{v}_V^{\frac12},
	\qquad && \forall\,v \in V,\\
    \label{Agm}
    \|v\|_{\infty}&\leq C\|v\|^{\frac12}\|v\|_{H^2(\Omega)}^{\frac12},
    \qquad &&\forall\,v\in H^2(\Omega).
\end{alignat}
Similar conclusions can be drawn for the corresponding vector- and matrix-valued spaces.

Let us now introduce the standard spaces for solenoidal vector fields (see, e.g., \cite{Sohr})
\begin{align*}
	\Hs &:=	\ov{\big\{ \bv \in C^\infty_c (\Omega)^2 : \, \div \bv =0\big\}}^{\LLL^2(\Omega)},
\quad \Vs := \ov{\big\{ \bv \in C^\infty_c (\Omega)^2 : \, \div \bv =0\big\}}^{\boldsymbol{H}^1(\Omega)}.
\end{align*}
They are endowed with natural norms indicated by $\norma{\cdot}_{\Hs}=\norma{\cdot}_{\LLL^2(\Omega)}$ and $\norma{\cdot}_{\Vs}=\norma{\cdot}_{\boldsymbol{H}^1(\Omega)}$, respectively.
The well-known Korn's inequality yields that
$$
\|\nabla \bu\| \leq \sqrt{2}\|D\bu\|,\quad  \forall\,\bu\in\Vs.
$$
This together with Poincar\'e's inequality implies that $\|D\cdot\|$ is an equivalent norm for $\Vs$. It is well known that $\boldsymbol{L}^2(\Omega)$ can be decomposed into $\Hs\oplus\boldsymbol{G}(\Omega)$, where $\boldsymbol{G}(\Omega):=\{\boldsymbol{g}\in\boldsymbol{L}^2(\Omega): \exists\, z\in H^1(\Omega),\ \boldsymbol{g}=\nabla z\}$. Then for any function $\boldsymbol{f} \in \boldsymbol{L}^2(\Omega)$, the Helmholtz--Weyl decomposition holds:
$$
\boldsymbol{f}=\boldsymbol{f}_{0}+\nabla z,\quad\text{where}\  \boldsymbol{f}_{0} \in \Hs,\ \nabla z \in \boldsymbol{G}(\Omega).\nonumber
$$
Consequently, we can define the Leray projection onto the space of divergence-free functions $\boldsymbol{P}:\boldsymbol{L}^2(\Omega)\to \Hs$ such that $\boldsymbol{P}(\boldsymbol{f})=\boldsymbol{f}_{0}$.
Introduce the Stokes operator $\boldsymbol{A}: D(\boldsymbol{A})= \boldsymbol{H}^2(\Omega) \cap \Vs \to\Hs$ such that $\boldsymbol{A}=\boldsymbol{P}(-\Delta)$ (see, e.g., \cite[Chapter III]{Sohr}).
We denote $\Ws := D(\boldsymbol{A})$, which can be equipped with the inner product $(\boldsymbol{A}\boldsymbol{u},\boldsymbol{A}\boldsymbol{v})$ and the equivalent norm $\|\boldsymbol{A}\boldsymbol{u}\|$.
Besides, we recall the following regularity result for the Stokes operator (see, e.g.,  \cite[Chapter III, Theorem 2.1.1]{Sohr}).
\begin{lemma}\label{sto}
Let $\Omega$ be a bounded domain of class $C^2$ in $\mathbb{R}^2$. For any $\boldsymbol{f} \in \Hs$,
there exists a unique pair $(\boldsymbol{u}, p)\in \Ws\times V_0$ such that $-\Delta \boldsymbol{u}+\nabla p=\boldsymbol{f}$ a.e. in $\Omega$, that is, $\boldsymbol{u}=\boldsymbol{A}^{-1}\boldsymbol{f}$. Moreover, it holds
\begin{align*}
&\|\boldsymbol{u}\|_{\boldsymbol{H}^2(\Omega)}+\|\nabla p\|\le C\|\boldsymbol{f}\|,\quad
\|p\|\le C \|\boldsymbol{f}\|^\frac12\|\nabla \boldsymbol{A}^{-1}\boldsymbol{f}\|^\frac12,
 \end{align*}
where $C>0$ depends on $\Omega$ but is independent of $\boldsymbol{f}$.
\end{lemma}

Finally, we set forth the following convention: the capital letters $C$ and $C_i$ indicate generic positive constants determined exclusively by the structure of the system. Its interpretation may vary from line to line and even within the same chain of computations. When a positive constant $\delta$ is part of the calculation, the corresponding symbol $C_\delta$ replaces the generic $C$, indicating a constant that also depends on $\delta$.

\subsection{Assumptions and results}
\label{SUBSEC:ASSMAINRES}

We start with some structural assumptions for the state system \eqref{eq:1}--\eqref{eq:7}.
Throughout this paper, we assume that
$$\eta\in \mathbb{R},\quad \varepsilon=1.$$
The precise value of the parameter $\varepsilon$ does not play an essential role in the subsequent analysis, because we do not consider here the sharp interface limit as $\varepsilon\to 0$.

Concerning the nonlinearities, we postulate the following.
\begin{enumerate}[label=$ \boldsymbol{(\mathrm{A \arabic*})}$, ref =$\boldsymbol {({\mathrm{A \arabic*}})}$]
\item \label{ass:1:visco}
		The viscosity function $\nu:\erre \to \erre^+$ is twice differentiable and
\begin{align*}
			\exists \,\nu_* >0 : \quad \nu(s)\geq \nu_*,  \quad  \forall\, s \in \erre.
\end{align*}
\item \label{ass:2:mob}
		The mobility function $m:\erre\to\erre^+$ is three times differentiable  and
\begin{align*}
			\exists \,m_*>0 : \quad m(s)\geq m_* , \quad  \forall\,  s \in \erre.
\end{align*}
\item \label{ass:3:pot}
The potential function $F:\erre\to\erre$ is five times differentiable. It can be decomposed as the sum of a convex function and its quadratic perturbation, that is, $F= F_1+F_2$ with $F_1$ being convex (so that $F_1''\geq 0$) and $F_2'$ being Lipschitz continuous.
Moreover, it holds that
\begin{alignat}{2}
    \label{growth:F:0}
	& \lim_{|s| \to +\infty} \frac{F_1(s)}{s^2} = + \infty,
\end{alignat}
and
\begin{alignat}{2}
           \label{growth:F:1-b}
     & \exists \, \gamma_1\in (0,1),\ \gamma_2 >0 \ \ \text{such that}\ \ sF'_1(s)\geq (2+\gamma_1) F_1(s)-\gamma_2,\quad \forall\, s\in\erre,
        \\
        \label{growth:F:2}
     & \exists \, c_F >0 \ \ \text{such that}\ \ |{F''_2(s)}|\le
        c_{F},\quad \forall\, s\in\erre.
    \end{alignat}
Without loss of generality, we assume $F_1(0)=F_1'(0)=0$.
\end{enumerate}

\begin{remark}\label{rem:bd}\rm
With the homogeneous Neumann boundary condition for $\varphi$, we find that $\dn \Delta \varphi =0$ is equivalent to $\dn \omega=0$ on $\Sigma$.
Next, assumption \ref{ass:2:mob} implies that the no-flux boundary condition $m(\ph)\dn \mu =0$ is equivalent to $\dn \mu =0$ on $\Sigma$. However, due to the nonlinear structure of $\mu$, condition $\dn \mu =0$ does not imply $\dn \Delta^2 \varphi =0$ on $\Sigma$.
\end{remark}

\begin{remark}\label{rem:F}\rm
The coercivity condition \eqref{growth:F:0} together with the growth assumption \eqref{growth:F:2} implies that $F$ is bounded from below on $\erre$. The regularity assumptions in \ref{ass:3:pot} are slightly overabundant to obtain only the weak and strong well-posedness of the state system (cf. \cite{WY}), but they are necessary for mathematical analysis of the associated control problem. On the other hand, \ref{ass:3:pot} is satisfied by a wide class of smooth potentials with either polynomial growth (of order greater than two) or exponential growth, in particular, by the classical regular potential \eqref{double-well}.
\end{remark}

To study the control problem \textbf{(CP)}, some additional assumptions are needed.
\begin{enumerate}[label=$ \boldsymbol{(\mathrm{A \arabic*})}$, ref =$\boldsymbol {(\mathrm{A \arabic*})}$, start=4]
\item \label{ass:4:uad}
The set of {\it admissible controls} $\Uad$ is given by
\begin{align*}
		\Uad = \big\{
	\bu \in \LLL^\infty(Q) \,\,:\,\, u_{{\rm min}, i}\leq u_{ i}	\leq u_{{\rm max}, i} \quad \aeQ, \quad i=1,2
	\big\},
\end{align*}
where $\umin=(u_{\min, 1},u_{\min, 2})$ and $\umax=(u_{\max, 1},u_{\max, 2})$ are certain given functions in $\LLL^\infty(Q)$ with $ u_{\min,i} \leq u_{\max,i}$ almost everywhere in $Q$, $i=1,2$.
Then we define the following open ball that contains $\Uad$:
\begin{align}
	\UR :=\{ \bu \in \LLL^\infty(Q): \,\, \norma{\bu}_{\LLL^\infty(Q)} <R \},
\label{UR}
\end{align}
with
$$
R:= \max_{i=1,2}\big\{\norma{u_{\min,i}}_{L^\infty(Q)},\, \norma{u_{\max,i}}_{L^\infty(Q)}\big\} +1.
$$
\item \label{ass:6:cost}
The real numbers $\alphaphQ,\alphaphO$ and $\alphacon$ in \eqref{def:cost} are nonnegative and not all zeros (to avoid a trivial problem). The target functions satisfy $\ph_Q \in L^2 (Q)$, $\ph_\Omega \in H$.
\end{enumerate}

Our first result concerns the strong well-posedness of the state system \Sys\ in the two-dimensional setting.
\begin{theorem}[Well-posedness]
\label{THM:RES}
Suppose that $T>0$, $\eta\in \mathbb{R}$, $\varepsilon=1$, $\bu \in \L2 {\HHH}$, the assumptions \ref{ass:1:visco}--\ref{ass:3:pot} are satisfied, and the initial data satisfy
$(\bv_0,\ph_0)\in \Vs \times (\Hx5\cap \Hqn)$.

({i}) Problem \eqref{eq:1}--\eqref{eq:7} admits a unique strong solution $\solst$ on $[0,T]$ such that
\begin{align*}
		\bv &  \in \L\infty {\Vs} \cap \L2 {\Ws}\cap \H1 {\Hs},
		\\
		p & \in \L2{V_0},
        \\
        \ph & \in
		\L\infty {\Hx5 \cap \Hqn} \cap \L2 {\Hx6},
		\\
		\dt \ph & \in \L\infty {\Vp} \cap \L2 {\Hdn},
		\\
		\mu & \in \L\infty {V} \cap \L2 {\Hdn},
		\\
		\omega & \in \L\infty {\Hx3 \cap \Hdn} \cap \L2 {\Hx4}.
\end{align*}
The solution $\solst$ satisfies the system \eqref{eq:1}--\eqref{eq:5} almost everywhere in $Q$, the boundary conditions \eqref{eq:6} almost everywhere on $\Sigma$ and the initial conditions \eqref{eq:7} almost everywhere in $\Omega$. Moreover, it holds
\begin{align}
	\non
	&
	\norma{\bv}_{ \L\infty{\Vs} \cap \L2 {\Ws}\cap \H1 {\Hs}}
 	+	\norma{p}_{\L2 V}
 	\\ & \quad \non
	+\norma{\ph}_{\L\infty{\Hx5} \cap \L2 {\Hx6}\cap \W{1,\infty}{\Vp} \cap \H1 {\Hdn} }
 	\\ & \quad \non
	+ \max_{i=0,...,5} \norma{F^{(i)}(\ph)}_{L^\infty(Q)} +\norma{\mu}_{\L\infty {V} \cap \L2 {\Hdn}}
 	\\ & \quad
 	 	+\norma{\omega}_{\L\infty {\Hx3 } \cap \L2 {\Hx4}}
 	\leq K_1,
 \label{control:reg}
\end{align}
where the positive constant $K_1$ depends on $\Omega$, $T$, the coefficients $\eta$, $\nu_*$, $m_*$, and the norms $\|\bv_0\|_{\Vs}$, $\|\ph_0\|_{H^5(\Omega)}$, $\|\bu\|_{\L2 {\HHH}}$.

({ii}) Given two strong solutions  $\{(\bv_i,p_i,\ph_i,\mu_i,\omega_i)\}$, $i=1,2$, associated to the same initial data $(\bv_0,\ph_0)$ and two forcing terms $\bu_i \in \L2 {\HHH}$, $i=1,2$, we have the following continuous dependence estimate
\begin{align}	
		&  \non
		\norma{\bv_1-\bv_2}_{C^0([0,t];\Hs)\cap L^2(0,t;\Vs)}
		+ \norma{\ph_1-\ph_2}_{C^0([0,t];\Hx2) \cap L^2(0,t;\Hx5)}
		\\ & \quad \label{cd:est:1}
		\leq K_2\norma{\boldsymbol{P}(\bu_1-\bu_2)}_{\L2 {\boldsymbol{V}_\sigma^*}}, \quad
\forall\, t\in(0,T],
\end{align}
for a positive constant $K_2$ that depends on $\Omega$, $T$, the coefficients $\eta$, $\nu_*$, $m_*$, and the norms $\|\bv_0\|$, $\|\ph_0\|_{H^2(\Omega)}$, $\|\bu_i\|_{\L2 {\boldsymbol{V}_\sigma^*}}$, $i=1,2$.
\end{theorem}

\begin{remark}\rm
Thanks to the Aubin--Lions--Simon compactness theorem (see \cite{Simon}) and the regularity properties obtained in Theorem \ref{THM:RES}-(i), we find $ \bv \in \C0{\Vs}$ and $\ph \in  \C0{\Hx4}$ so that the initial data can be attained.
\end{remark}
\begin{remark}\rm
\label{RMK:pressure}
The pressure $p$ can be recovered a posteriori, up to an additive constant, like, e.g., in  \cite[Chapter V, Section 1.5]{Bo21} such that $p\in W^{-1,\infty}(0,T; H_0)$ (see also \cite[Chapter IV, Lemma 1.4.1]{Sohr}).
Then from the regularity properties of the strong solution $(\bv, \ph, \mu)$ and the assumption \ref{ass:1:visco}, we have
\begin{alignat*}{2}
	&  \nabla p
	= - \dt \bv - (\bv \cdot \nabla ) \bv + 2 \mathrm{div}(\nu (\ph) D \bv)
      +\mu \nabla \ph + \bu\in \L2{\HHH},
\end{alignat*}
which implies $p \in \L2{V_0}$ (since $p$ has a zero mean, so the Poincar\'e--Wirtinger inequality applies).
\end{remark}

Using the expression of the chemical potentials $\mu$ and $\omega$, we can obtain the following continuous dependence estimate with respect to $\bu$.

\begin{corollary}[Continuous dependence for chemical potentials] \label{THM:CD:NEW}
	Suppose that the assumptions of Theorem \ref{THM:RES} are fulfilled. Then, it holds
\begin{align}	
		& \norma{\mu_1-\mu_2}_{L^2(0,t;V)}
		+ \norma{\omega_1-\omega_2}_{C^0([0,t]; H) \cap L^2(0,t;\Hx3)} \non\\
        & \quad
		\label{cd:est:2}
		\leq K_3\norma{\boldsymbol{P}(\bu_1-\bu_2)}_{\L2 {\boldsymbol{V}_\sigma^*}}, \quad \forall\, t\in (0,T],
\end{align}
for a positive constant $K_3$  that depends on $\Omega$, $T$, the coefficients $\eta$, $\nu_*$, $m_*$, and the norms $\|\bv_0\|$, $\|\ph_0\|_{H^2(\Omega)}$, $\|\bu_i\|_{\L2 {\boldsymbol{V}_\sigma^*}}$, $i=1,2$.
\end{corollary}

In control theory, the external variable $\bu$ occurring in \eqref{eq:1} is termed the {\it control}, and the corresponding solution $\solst$ to problem \Sys\ is termed the {\it state}. Then Theorem \ref{THM:RES} and Corollary \ref{THM:CD:NEW} allow us to define the {\it control-to-state} mapping:
\begin{align*}
	\Sol: \ \UR &\to \YY,\\
          \bu & \mapsto \solst{,}
\end{align*}
where $\YY$ stands for the solution space emerging in Theorem \ref{THM:RES} and the open set $\UR$ is defined in \eqref{UR}.
\begin{remark}\rm
Due to our choice of the cost functional \eqref{def:cost}, it will be more convenient to interpret the above control-to-state mapping simply as
$$\Sol: \bu \mapsto \ph,$$
with $\ph$ being the third component of the unique strong solution $\solst$ to problem \Sys. As this also improves the clarity of the presentation, we shall proceed in this manner, occasionally without explicitly stating whether the operator $\Sol$ refers only to the phase-field variable $\ph$ or to all variables.
\end{remark}

Define the Banach space
\begin{align}
	\label{def:Z}
		\ZZ:=
 \C0 {\Hdn } \cap \L2 {\Hx5\cap \Hqn}\cap \H1 \Vp.
\end{align}
It follows from \eqref{cd:est:1} that the control-to-state operator $\Sol$, i.e., $\bu \mapsto \Sol(\bu)=\ph$,
is locally \Lip\ continuous from $\L2 {\Vs^*}$ into $\ZZ$.

In order to establish the differentiability of $\Sol$,
we shall study the following \textit{linearized system}. Let $\bu \in \UR$ be a control with the corresponding state $\solst$ determined by Theorem \ref{THM:RES}. For every given $\bh \in \L2 \HHH$, we consider
\begin{align}
	\non
	& \dt \bw
	+ (\bw \cdot \nabla ) \bv
	+ (\bv \cdot \nabla ) \bw
	- \div (2 \nu' (\ph) \psi D \bv)
	- \div (2 \nu (\ph) D \bw)
	+ \nabla q
	&&
	\\ 	\label{eq:lin:1}
	& \quad
	=
	\theta \nabla \ph
	+ \mu \nabla \psi
	+ \bh
	\qquad && \text{in $Q$},
	\\
	\label{eq:lin:2}
	& \div \bw =0
	\qquad && \text{in $Q$},
	\\
	\label{eq:lin:3}
	& \dt \psi
	+ \bw \cdot \nabla \ph
	+ \bv \cdot \nabla \psi
	- \div (m'(\ph)\psi\nabla \mu)
	- \div (m(\ph)\nabla \theta)
	=0
	\qquad && \text{in $Q$},
	\\
	\label{eq:lin:4}
	& \theta = - \Delta \w
	+ f''(\ph)\psi \omega
	+ f'(\ph) \w
	+ \eta \w
	\qquad && \text{in $Q$},
	\\
	\label{eq:lin:5}
	& \w =
	-\Delta \psi
	+ f'(\ph)\psi
	\qquad && \text{in $Q$},
		\\
	\label{eq:lin:6}
	& \bw = {\boldsymbol 0},
	\quad
	\dn \psi = \dn \Delta \psi = \dn \theta =0
	\qquad && \text{on $\Sigma$},
	\\
	\label{eq:lin:7}
	& \bw(0)=\0,
	\quad
	\psi(0)=0
	\qquad && \text{in $\Omega$}.
\end{align}
%
As observed in Remark \ref{rem:bd}, the homogeneous Neumann boundary condition $\dn \Delta \psi =0$ is equivalent to $\dn \w=0$ on $\Sigma$.
Well-posedness of the linear evolution problem \eqref{eq:lin:1}--\eqref{eq:lin:7} on $[0,T]$ will be demonstrated in Proposition \ref{THM:LIN}.
With the aid of this property, we can show that the control-to-state operator $\Sol$ is Fr\'{e}chet differentiable between suitable Banach spaces, and establish a connection between its derivative and the solution to the linearized system \eqref{eq:lin:1}--\eqref{eq:lin:7}. The result is summarized in the following theorem.

\begin{theorem}[Differentiability of the control-to-state operator]
\label{THM:FRE}
Suppose that $T>0$, $\eta\in \mathbb{R}$, $\varepsilon=1$, $\bu \in \L2 {\HHH}$, the assumptions \ref{ass:1:visco}--\ref{ass:4:uad} are satisfied, and the initial data satisfy
$(\bv_0,\ph_0)\in \Vs \times (\Hx5\cap \Hqn)$.
The control-to-state operator $\Sol$ is \Fre\ differentiable in $\UR$ as a mapping from $\UR$ into $\ZZ$. Moreover, for any $\bu \in \UR$, it holds that
	$D\Sol (\bu) \in {\cal L}(\UR, \ZZ)$ and $$D\Sol (\bu)[\bh] = \psi,$$ where $\psi$ is the third component of $\sollin$ that represents the unique solution to the linearized system \eqref{eq:lin:1}--\eqref{eq:lin:7} associated with the control $\bu$ and the variation $\bh$ as in Proposition \ref{THM:LIN}.
\end{theorem}

Let us now delve into results about the optimal control problem \textbf{(CP)}.
The first result shows that the minimization problem \textbf{(CP)} under investigation admits at least one solution, which is usually termed an {\it optimal control}.

\begin{theorem}[Existence of an optimal control]
	\label{THM:EX:CONTROL}
Suppose that $T>0$, $\eta\in \mathbb{R}$, $\varepsilon=1$, the assumptions \ref{ass:1:visco}--\ref{ass:6:cost} are satisfied, and the initial data satisfy $(\bv_0,\ph_0)\in \Vs \times (\Hx5\cap \Hqn)$.
Then, the optimal control problem \textbf{(CP)} admits at least one (globally) optimal solution. Namely, there exists an admissible control variable $\opt\in \Uad$ such that
\begin{align*}
	 {\cal J}(\opt,\Sol(\opt)) \leq {\cal J}(\bu,\Sol(\bu)),
  \quad \forall\, \bu\in\Uad.
\end{align*}
\end{theorem}

Due to the highly nonlinear nature of the state system \eqref{eq:1}--\eqref{eq:7}, uniqueness of the optimal control $\opt$ is not expected in general. However, by exploiting the differentiability property of the control-to-state operator $\Sol$ (recall Theorem \ref{THM:FRE}), we can derive first-order necessary optimality conditions associated with the optimal control problem.

\begin{theorem}[First-order necessary optimality condition]	
\label{THM:VAR:INEQ:PREL}
Suppose that the assumptions in Theorem \ref{THM:EX:CONTROL} are satisfied. 	Let $\opt\in \Uad$ be an optimal control with the corresponding state $(\widetilde{\bv},\widetilde{p},\bph, \widetilde{\mu},\widetilde{\omega})$ given by Theorem \ref{THM:RES}.
Then the following variational inequality holds
\begin{align}
		& \alphaphQ \intQ (\bph- \ph_Q) \psi \,\mathrm{d}x\mathrm{d}t
		+ \alphaphO	\iO (\bph(T)- \ph_\Omega) \psi(T)\,\mathrm{d}x \notag \\
&\quad
		+ \alphacon \intQ \opt \cdot (\bu-\opt)\,\mathrm{d}x\mathrm{d}t
		\geq 0,
		\quad \forall\,\bu \in  \Uad,
\label{var:ineq:prel}
\end{align}
where $\psi$ is the third component of the associated variables $\sollin$ satisfying the linearized system  \eqref{eq:lin:1}--\eqref{eq:lin:7} with
$(\bv, p, \ph, \mu, \omega)=(\widetilde{\bv},\widetilde{p},\bph, \widetilde{\mu},\widetilde{\omega})$ and $\bh = \bu-\opt$.
\end{theorem}

The variational inequality \eqref{var:ineq:prel} is not satisfactory, since it requires solving the linearized system \eqref{eq:lin:1}--\eqref{eq:lin:7}  with $\bh = \bu-\opt$ infinitely many times as $\bu$
varies in $\Uad$. In optimal control theory, this limitation can be overcome by working with the so-called {\it adjoint system}, whose solutions play the role of Lagrangian multipliers. For simplicity, we shall present the adjoint system in its strong form, even if the problem will be solved only in a suitable weak sense (see Section \ref{SEC:CONTROL}).

Let $\bu \in \UR$ with $\solst$ be the corresponding state given by Theorem \ref{THM:RES}. Then, the adjoint system consists in looking for the tuple $\soladj$ solving (in a weak sense) the following backward system:
%
\begin{alignat}{2}
	\label{eq:ad:1}
	& -\dt \bva
	- \div (2 \nu(\ph)D\bva)
	- (\bv\cdot \nabla)\bva
	+ (\bva \cdot \nabla^\top)\bv
	+ \pha \nabla \ph
	- \nabla p^a
	=
	0
	\qquad && \text{in $Q$},
	\\
	\label{eq:ad:2}
	& \div \bva =0
	\qquad && \text{in $Q$},
	\\
	\non
	& -\dt \pha
	- \Delta \oma
	+ 2 \nu'(\ph)D\bv : \nabla \bva
	- \nabla \pha\cdot \bv
	+ \nabla \mu \cdot \bva
	+m'(\ph)\nabla \mu \cdot \nabla \pha
	&&
	\\
	& \quad
	\label{eq:ad:3}
	+ f''(\ph)\omega \mua
	+f'(\ph)\oma
	=\alphaphQ(\ph-\ph_Q)
	\qquad && \text{in $Q$},
	\\
	\label{eq:ad:4}
	&
	\mua =
	- \div (m(\ph) \nabla \pha)
	- \nabla  \ph \cdot \bva
	\qquad && \text{in $Q$},
	\\
	\label{eq:ad:5}
	& \oma
	=
	- \Delta \mua
	+f'(\ph)\mua
	+\eta \mua
	\qquad && \text{in $Q$},
		\\
	\label{eq:ad:6}
	& \bva = {\boldsymbol 0},
	\quad
	\dn \pha = \dn \mua = \dn \oma = 0
	\qquad && \text{on $\Sigma$},
	\\
	\label{eq:ad:7}
	& \bva(T)={\boldsymbol 0},
	\quad
	\pha(T)=\alphaphO (\ph(T)-\ph_\Omega)
	\qquad && \text{in $\Omega$}.
\end{alignat}
\Accorpa\Adj {eq:ad:1} {eq:ad:7}
The adjoint system \Adj\ can be derived by the formal Lagrange method, see, e.g., \cite{Tr10}. In fact, it can also be deduced directly from the calculations in the proof of Theorem \ref{THM:VAR:INEQ:FINAL} stated below. In order to simplify the presentation and keep focus on the optimal control problem, we assume constant mobility in the subsequent analysis, that is,
 $$m(\ph)\equiv 1.$$
Further comments on this aspect can be found in Section \ref{SUBSEC:FOC}.

The well-posedness of the adjoint system \Adj\ will be established in Proposition \ref{THM:ADJ}. With the aid of the solution to problem \Adj, we can simplify the variational inequality \eqref{var:ineq:prel} and establish the following first-order necessary optimality condition for the control problem \textbf{(CP)}.

\begin{theorem}[Optimality condition via adjoint states]	
\label{THM:VAR:INEQ:FINAL}
Suppose that the assumptions in Theorem \ref{THM:EX:CONTROL} are satisfied and, in addition, $m(\ph)\equiv 1$. Let $\opt$ be an optimal control with the corresponding state $(\widetilde{\bv}, \widetilde{p}, \bph,\widetilde{\mu}, \widetilde{\omega})$ given by Theorem \ref{THM:RES} and the adjoint state $(\bva,p^a,\pha,\mua,\oma)$ given by Proposition \ref{THM:ADJ}.
Then the following variational inequality holds
\begin{align}
		\label{var:ineq}
		\intQ (\alphacon \opt  + \bva )\cdot ( \bu-\opt)\,\mathrm{d}x\mathrm{d}t \geq 0,
		\quad \forall\, \bu \in  \Uad.
\end{align}
Furthermore, whenever $\alphacon$ is strictly positive, the optimal control can be determined by
$$
\opt = {\PPP}_{\Uad}(-{\alphacon}^{-1} \bva),
$$
where $\PPP_{\Uad}$ denotes the $\LLL^2$-orthogonal projection onto the convex set $\Uad$.
\end{theorem}

\begin{remark}\rm
\comm{Comparing \eqref{defE} and \eqref{energy-Pen}, we find that the free energies considered in this study and in \cite{ClimentEzquerraGuillenGonzalez2016} differ only in lower-order terms. Thus, our results on the well-posedness of problem \eqref{eq:1}--\eqref{eq:7} (Theorem \ref{THM:RES}, Corollary \ref{THM:CD:NEW}), the differentiability of the control-to-state operator (Theorem \ref{THM:FRE}), the existence of an optimal control (Theorem \ref{THM:EX:CONTROL}) as well as the first-order necessary optimality conditions (Theorem \ref{THM:VAR:INEQ:PREL}, Theorem \ref{THM:VAR:INEQ:FINAL}) can be correspondingly extended to the Cahn--Hilliard--Navier--Stokes vesicle-fluid interaction model analyzed in \cite{ClimentEzquerraGuillenGonzalez2016}. We leave the details to interested readers.}
\end{remark}

\section{The Control-to-State Operator and Its Properties}
\label{SEC:DIFF}
\setcounter{equation}{0}

Let us begin with a brief overview of the analytical strategy of this section. Recall that our goal is to identify appropriate first-order necessary conditions for optimality. Formally, this requires differentiating the cost functional $\cal J$ given by \eqref{def:cost}. Since $\cal J$ is quadratic, it suffices to analyze the properties of the solution operator $\Sol$, which constitutes the main source of nonlinearity.
The first step is to show that the solution operator $\Sol$ is Fr\'{e}chet differentiable in a suitable functional framework, so that the chain rule in Banach spaces can be applied to derive a preliminary characterization of the necessary optimality conditions in the form of a variational inequality, as specified in Theorem \ref{THM:VAR:INEQ:PREL}.
Subsequently, in Section \ref{SEC:CONTROL}, we shall simplify this condition using solutions to the adjoint system, as stated in Theorem \ref{THM:VAR:INEQ:FINAL}.
Based on these considerations, in this section, we  establish some key properties of the control-to-state operator $\Sol$.

\subsection{Well-posedness and stability of the state system}
We begin with the proof of Theorem \ref{THM:RES} on the strong well-posedness of problem \eqref{eq:1}--\eqref{eq:7}.
\begin{proof}[Proof of Theorem \ref{THM:RES}]
In view of the proofs for \cite[Theorems 2.5, 2.7, 2.8]{WY}, we only sketch the proof by pointing out the main differences and necessary modifications in the two-dimensional case.

(1) \textbf{Existence of strong solutions}.
The proof of existence can be done by means of a Faedo--Galerkin scheme similar to that in \cite[Section 3]{WY}. In the following, we just perform \textit{a priori} estimates for sufficiently smooth solutions along with the additional forcing term $\bu \in \L2 \HHH$.

\textit{Mass conservation}. Testing \eqref{eq:3} by $1$, using integration by parts over $\Omega$ and then integration on $[0,t]\subset [0,T]$, we easily find the mass conservation law
\begin{align}
\int_\Omega \ph(t)\,\mathrm{d}x=\int_\Omega \ph_0\,\mathrm{d}x,\quad \forall\, t\in [0,T].\label{app-e0}
\end{align}

\textit{Basic energy estimate}.  Following a similar argument for \cite[Lemma 3.1]{WY}, testing
\eqref{eq:1} by $\bv$ and \eqref{eq:3} by $\mu$, respectively, adding the
resultants together, integrating over $\Omega$ and then on $[0,t]\subset[0,T]$, we can derive the following energy equality
\begin{align}
&\frac{1}{2}\|\bv(t)\|^2+{\cal E}(\ph(t))
 +\int_0^t \int_\Omega \big(2\nu(\ph(\tau))|D\bv(\tau)|^2 +m(\ph(\tau)) |\nabla \mu(\tau)|^2\big)\,\mathrm{d}x \mathrm{d}\tau
\notag\\
&\quad = \frac{1}{2}\|\bv_0\|^2+{\cal E}(\ph_0) +
\int_0^t \int_\Omega \bu(\tau)\cdot \bv(\tau)\,\mathrm{d}x\mathrm{d}\tau,
\quad \forall\, t\in [0,T].
\label{app-e1}
\end{align}
Here, the free energy ${\cal E}$ is defined as in \eqref{defE}.
By the Sobolev embedding theorem in two dimensions and \ref{ass:3:pot}, the initial free energy ${\cal E}(\ph_0)$ can be bounded by
a positive constant $C$ depending on $\Omega$, $\eta$, and $\|\ph_0\|_{H^2(\Omega)}$.

Next, we show that ${\cal E}$ is bounded from below by a uniform constant. The case $\eta\geq 0$ is trivial since $F(s)$ is bounded from below for $s\in \mathbb{R}$ (see Remark \ref{rem:F}). Thus, it remains to treat the case $\eta<0$. We infer from the assumption \eqref{growth:F:1-b} that
$$
sF_1'(s)-2F_1(s)\geq \gamma_1F_1(s)-\gamma_2,\quad \forall\, s\in \mathbb{R}.
$$
Besides, the growth condition \eqref{growth:F:2} yields
$$
|sF_2'(s)-2F_2(s)|\leq C_1s^2+C_2,\quad \forall\, s\in \mathbb{R},
$$
for some positive constants $C_1, C_2$. Hence, from \ref{ass:3:pot}, the decomposition $f(\ph)=F_1'(\ph)+F_2'(\ph)$ and Young's inequality, we get
\begin{align*}
& \eta \int_\Omega \Big(\frac{1}{2}|\nabla \ph|^2 +F( \ph)\Big) \,\mathrm{d}x\\
&\quad = \frac{\eta}{2} \int_\Omega  \ph (-\Delta \ph+f( \ph))\,\mathrm{d}x
- \frac{\eta}{2} \int_\Omega ( \ph F_1'( \ph)-2F_1( \ph))\,\mathrm{d}x
- \frac{\eta}{2} \int_\Omega ( \ph F_2'( \ph)-2F_2( \ph))\,\mathrm{d}x
\notag \\
&\quad \geq -\frac{1}{4} \int_\Omega (-\Delta  \ph +f( \ph))^2\,\mathrm{d} x
-\frac{\eta^2}{4}\int_\Omega  \ph^2\,\mathrm{d}x
-\frac{\eta}{2}\int_\Omega (\gamma_1F_1(\ph)-\gamma_2)\,\mathrm{d}x \notag\\
&\qquad +\frac{\eta}{2}\int_\Omega (C_1\varphi^2+C_2)\,\mathrm{d}x \notag\\
&\quad \geq  -\frac{1}{4} \int_\Omega (-\Delta \ph +f(\ph))^2\,\mathrm{d} x -\frac{\eta\gamma_1}{4}\int_\Omega F_1(\ph)\,\mathrm{d}x
-C_3,
\end{align*}
where $C_3>0$ depends on $|\Omega|$, $\eta$, $\gamma_1$, $\gamma_2$, $C_1$, $C_2$, but is independent of $\ph$. In summary, for $\eta\in \mathbb{R}$, it holds (recall \eqref{eq:5}, \eqref{defE})
\begin{align}
& {\cal E}(\ph)\geq \frac{1}{4} \|\omega\|^2 -C_4,
\label{app-e3}
\end{align}
where $C_4>0$ may depend on $|\Omega|$, $\eta$, $\gamma_1$, $\gamma_2$, $C_1$, $C_2$, but it is independent of $\ph$.

For the last term on the right-hand side of \eqref{app-e1}, using Korn's inequality, Poincar\'e's inequality and Young's inequality, we deduce from \ref{ass:1:visco} that
\begin{align}
\left|\int_0^t \int_\Omega \bu(\tau)\cdot \bv(\tau)\,\mathrm{d}x\mathrm{d}\tau\right|
&\leq \nu_*\int_0^t \int_\Omega |D\bv(\tau)|^2\,\mathrm{d}x \mathrm{d}\tau
+ \frac{C}{\nu_*} \int_0^t \int_\Omega |\bu(\tau)|^2\,\mathrm{d}x\mathrm{d}\tau,
\label{app-e3b}
\end{align}
where $C>0$ depends only on $\Omega$.

Combining \eqref{app-e1}, \eqref{app-e3} and \eqref{app-e3b}, we can conclude
\begin{align}
\|\bv(t)\|^2+ \|\omega(t)\|^2 + \int_0^t \int_\Omega \big(\nu_*|D\bv(\tau)|^2 +m_* |\nabla \mu(\tau)|^2\big)\,\mathrm{d}x \mathrm{d}\tau\leq C, \quad \forall\, t\in [0,T],
\label{app-energy-low}
\end{align}
for some $C>0$ depending on $\|\bv_0\|$, $\|\ph_0\|_{H^2(\Omega)}$, $\|\bu\|_{L^2(0,T;\boldsymbol{H})}$, $\Omega$ and structural constants of the system.

To obtain an $H^2$-estimate for $\ph$, we consider the following Neumann problem
\begin{equation}
\begin{cases}
-\Delta \ph + \ph  = \omega + \ph - f(\ph)\qquad \text{in}\ \Omega,\\
\partial_\mathbf{n} \ph=0,\qquad\qquad \qquad \qquad\hspace{.5cm}\text{on}\ \Gamma.
\end{cases}
\label{ellipA}
\end{equation}
Multiplying the equation in \eqref{ellipA} by $\ph$, integrating over $\Omega$, we infer from \eqref{growth:F:0}--\eqref{growth:F:2} that
\begin{align}
\|\ph\|_{H^1(\Omega)}^2
&=  \int_\Omega \omega \ph\,\mathrm{d}x +\|\ph\|^2
- \int_\Omega \ph F_1'(\ph)\,\mathrm{d}x
- \int_\Omega \ph F_2'(\ph)\,\mathrm{d}x
\notag \\
&\leq  \|\omega\|^2 + \frac54\|\ph\|^2 - (2+\gamma_1)\int_\Omega F_1(\ph)\,\mathrm{d}x + \gamma_2|\Omega| + C_5\|\ph\|^2+C_6
\notag\\
&\leq \|\omega\|^2 +C_7,
\label{app-energy-H1}
\end{align}
where $C_7>0$ depends on $|\Omega|$, $\eta$, $\gamma_1$, $\gamma_2$, $C_5$, $C_6$, but is independent of $\ph$. On the other hand, testing the equation in \eqref{ellipA} by $F_1'(\ph)$, we easily find
\begin{align}
\underbrace{\int_\Omega F_1''(\ph)|\nabla \ph|^2\,\mathrm{d}x}_{\geq 0} \,+\, \|F_1'(\ph)\|^2
&  = \int_\Omega (\omega-F_2'(\ph)) F_1'(\ph)\,\mathrm{d}x
\notag \\
&\leq \frac12 \|F_1'(\ph)\|^2 + \|\omega\|^2 + \|F_2'(\ph)\|^2
\notag \\
&\leq \frac12 \|F_1'(\ph)\|^2 + \|\omega\|^2 + C(\|\ph\|^2+1).
\notag
\end{align}
Applying the elliptic estimate for the Neumann problem \eqref{ellipA}, we infer from the above estimate that
\begin{align}
\|\ph\|_{H^2(\Omega)} &\leq C(\|\omega-f(\ph)\|+\|\ph\|)
\notag \\
&\leq C(\|\omega\|+\|F_1'(\ph)\|+ \|F_2'(\ph)\| +\|\ph\|)
\notag\\
& \leq C(\|\omega\|+\|\ph\|+1).
\label{app-e4}
\end{align}
This combined with \eqref{app-energy-low} and \eqref{app-energy-H1} yields
\begin{align}
\|\ph(t)\|_{H^2(\Omega)}\leq C,\quad \forall\, t\in [0,T].
\label{app-energy-H2}
\end{align}
Noticing that
\begin{align}
\left|\int_\Omega \mu\,\mathrm{d}x\right|\leq (\|F_1'(\ph)\|+\|F_2'(\ph)\|)\|\omega\|+ |\eta||\Omega|^\frac12\|\omega\|\in L^\infty(0,T),
\label{mu-mean}
\end{align}
then, arguing as in \cite[Section 3]{WY}, we can apply \eqref{app-energy-low}, \eqref{app-energy-H2} to obtain
\begin{align}
\int_0^t \|\mu(\tau)\|_{H^1(\Omega)}^2\,\mathrm{d}\tau+ \int_0^t \|\omega(\tau)\|_{H^3(\Omega)}^2\,\mathrm{d}\tau
+\int_0^t \|\ph(\tau)\|_{H^5(\Omega)}^2\,\mathrm{d}\tau\leq C,\quad  \forall\, t\in [0,T].\label{app-e7}
\end{align}

\textit{Higher-order energy estimate}. Using the Sobolev embedding theorem in two dimensions and the lower-order estimate \eqref{app-energy-low}, we get
\begin{align}
|(\bv\cdot\nabla \bv, \boldsymbol{A}\bv)|\leq \|\bv\|_{4}\|\nabla\bv\|_{4} \|\boldsymbol{A}\bv\| \leq  C\|\bv\|^\frac12 \|\nabla \bv\| \|\boldsymbol{A}\bv\|^\frac{3}{2}  \leq \dfrac{\nu_*}{12}\|\boldsymbol{A}\bv\|^2 + C\|\nabla \bv\|^4.
\label{term 1}
 \end{align}
Keeping \eqref{term 1} in mind, by a careful examination of the argument in \cite[Section 4]{WY}, we can derive the following differential inequality
\begin{align}
&\frac{\mathrm{d}}{\mathrm{d}t} \left(\frac12 \|\nabla \bv\|^2+\frac{\lambda}{2}\int_\Omega m(\ph)|\nabla \mu|^2\,\mathrm{d}x +\lambda \big(\bv\cdot\nabla \ph,\mu\big)\right)\notag \\
&\qquad +\frac{\nu_*}{2}\|\boldsymbol{A}\bv\|^2 + \frac{\nu_*}{8C_*}\|\partial_t\bv\|^2 +\frac{3\lambda}{4} \|\Delta\partial_t\varphi\|^2\notag\\
&\quad \leq C\|\nabla \bv\|^4 + C\lambda \|\nabla \mu\|^4 + C(1+\lambda)(\|\nabla \bv\|^2 +\|\nabla \mu\|^2) +C\|\bu\|^2,
\label{duu-mu}
\end{align}
where $C_*, C, \lambda$ are positive constants depending on $\|\bv_0\|$, $\|\ph_0\|_{H^2(\Omega)}$, $\|\bu\|_{L^2(0,T;\boldsymbol{H})}$, $\Omega$ and structural constants of the system. Here, we just point out that to derive \eqref{duu-mu}, different types of interpolation in two dimensions are used compared to the original argument in \cite{WY} for three dimensions. Set
\begin{align*}
\Lambda(t)
& = \frac12 \|\nabla \bv(t)\|^2+\frac{\lambda}{2}\int_\Omega m(\ph(t))|\nabla \mu(t)|^2\,\mathrm{d}x
+ \lambda (\bv(t)\cdot\nabla \ph(t),\mu(t)).
\end{align*}
Thanks to \eqref{app-energy-H2}, there exists some positive constant $m^*$ satisfying
$m^*>m_*$ and $\|m(\ph)\|_{L^\infty(\Omega)}\leq m^*$.
As a consequence, it holds
\begin{align}
&\frac14\|\nabla \bv(t)\|^2+ \frac{\lambda m_*}{4}\|\nabla \mu(t)\|^2\leq \Lambda(t)\leq \|\nabla \bv(t)\|^2+ \lambda m^*\|\nabla \mu(t)\|^2.\label{LL}
\end{align}
Thus, we can deduce from \eqref{duu-mu} and \eqref{LL} that
\begin{align}
& \frac{\mathrm{d}}{\mathrm{d}t}\Lambda(t)  +\frac{\nu_*}{2}\|\boldsymbol{A}\bv\|^2 + \frac{\nu_*}{8C_*}\|\partial_t\bv\|^2 +\frac{3\lambda}{4} \|\Delta\partial_t\varphi\|^2
 \leq \widetilde{C}\Lambda(t)\big(\Lambda(t)+1\big) + C\|\bu\|^2,
\label{La}
\end{align}
where the positive constant $\widetilde{C}$ depends on $\|\bv_0\|$, $\|\ph_0\|_{H^2(\Omega)}$, $\|\bu\|_{L^2(0,T;\boldsymbol{H})}$, $\Omega$ and structural constants of the system.

The assumptions on the initial data yield $\Lambda(0)\leq C$, where $C>0$ depends on $\|\bv_0\|_{\boldsymbol{V}_\sigma}$, $\|\ph_0\|_{H^5(\Omega)}$ and $\Omega$. In addition, noticing \eqref{app-energy-low} and the following fact
\begin{align*}
|(\bv\cdot\nabla \ph,\mu)| &=|(\bv\cdot\nabla \ph,\mu-\overline{\mu})| \leq \|\bv\|_{3}\|\nabla \ph\|_{3}\|\mu-\overline{\mu}\|_{3}\\
&\leq C\|\nabla \ph\|_{3}\|\nabla \bv\|\|\nabla \mu\| \leq \frac{C}{4m_*}\|\nabla \bv\|^2+\frac{m_*}{4}\|\nabla \mu\|^2,
\end{align*}
we find $ \Lambda(t)\in L^1(0,T)$.
From these observations and the fact $\bu\in L^2(0,T;\boldsymbol{H})$, we can apply Gronwall's lemma to \eqref{La} and conclude
$$
\Lambda(t)+\int_0^{T} \left(\|\boldsymbol{A}\bv(t)\|^2 + \|\partial_t\bv(t)\|^2 + \|\Delta\partial_t\varphi(t)\|^2\right)\mathrm{d}t\leq C,
\quad \forall\, t\in [0,T],
$$
which together with \eqref{app-e0}, \eqref{app-energy-low}, \eqref{app-e4}, \eqref{mu-mean},  \eqref{LL} yields the following estimates
\begin{align}
& \|\bv\|_{L^\infty(0,T;\boldsymbol{V}_\sigma)}+ \|\bv\|_{L^2(0,T;\boldsymbol{W}_\sigma)} +
\|\partial_t\bv\|_{L^2(0,T;\boldsymbol{H}_\sigma)}\leq C,
\label{app-e8}
\\
& \|\mu\|_{L^\infty(0,T;H^1(\Omega))}\leq C,\label{app-e9}
\\
& \|\partial_t \ph\|_{L^2(0,T;H^2(\Omega))}\leq C,\label{app-e10}
\end{align}
where $C>0$ depends on $\|\bv_0\|_{\boldsymbol{V}_\sigma}$, $\|\ph_0\|_{H^5(\Omega)}$, $\|\bu\|_{L^2(0,T;\boldsymbol{H})}$, $\Omega$ and structural constants of the system. Furthermore, by the same argument as in \cite[Section 4]{WY}, we can further obtain
\begin{align}
&\|\ph\|_{L^\infty(0,T;H^5(\Omega))} + \|\omega\|_{L^\infty(0,T;H^3(\Omega))}\leq C,\label{app-e11}
\\
&\|\ph\|_{L^2(0,T;H^6(\Omega))}
+ \|\omega\|_{L^2(0,T;H^4(\Omega))}
+ \|\mu\|_{L^2(0,T;H^2(\Omega))}\leq C.
\label{app-e12}
\end{align}

With the above estimates, we can work with a suitable Galerkin approximation and apply the compactness method to establish the existence of a global strong solution to problem \eqref{eq:1}--\eqref{eq:7} on $[0,T]$.
\medskip

(2) \textbf{Continuous dependence and uniqueness}. The uniqueness of strong solutions to problem \eqref{eq:1}--\eqref{eq:7} is a direct consequence of \eqref{cd:est:1}.
In order to show this stability estimate with respect to perturbations of the external forcing term $\bu$, we exploit the proof of the weak-strong uniqueness result \cite[Theorem 2.5]{WY}, which is based on the relative energy method. To handle the additional forcing term, we apply the following fact
\begin{align*}
\left|\int_0^t\int_\Omega (\bu_1-\bu_2)\cdot(\bv_1-\bv_2)\,\mathrm{d}x\mathrm{d}\tau\right|
&\leq \int_0^t\|\boldsymbol{P}(\bu_1-\bu_2)\|_{\boldsymbol{V}_\sigma^*}
\|\bv_1-\bv_2\|_{\boldsymbol{V}_\sigma}\,\mathrm{d}\tau
\\
&\leq \delta \int_0^t \|D(\bv_1-\bv_2)\|^2\,\mathrm{d}\tau
+C_\delta\int_0^t\|\boldsymbol{P}(\bu_1-\bu_2)\|_{\boldsymbol{V}_\sigma^*}^2 \,\mathrm{d}\tau,
\end{align*}
for any $t\in [0,T]$. Here, we have used Poincar\'{e}'s inequality, Korn's inequality, Young's inequality, and chosen $\delta\in (0,1)$ as a suitably small constant. Then, by the Sobolev embedding theorem in two dimensions, we can derive the continuous dependence estimate \eqref{cd:est:1} following the argument in \cite{WY} with minor modifications.
The details are omitted here.

The proof of Theorem \ref{THM:RES} is complete.
\end{proof}

In what follows, we sketch the proof of Corollary \ref{THM:CD:NEW}.
\begin{proof}[Proof of Corollary \ref{THM:CD:NEW}]
The proof of the continuous dependence estimate \eqref{cd:est:2} for $(\mu, \omega)$ easily follows from a comparison argument in the difference of equations \eqref{eq:4} and \eqref{eq:5} along with the continuous dependence estimate \eqref{cd:est:1} for $\varphi$ and a technical lemma \cite[Lemma 3.2]{WY} (with a slight extension under assumption \ref{ass:3:pot}). Thus, we omit the details.
\end{proof}

\subsection{Differentiability of the control-to-state operator}
Let us now move to the Fr\'{e}chet differentiability of the solution operator $\Sol$ associated with the state system \Sys. This property will be essential for investigating the optimal control problem \textbf{(CP)}. To establish the forthcoming results, we rely on estimates derived from the higher-order regularities in Theorem \ref{THM:RES}.

\begin{proposition}
\label{THM:LIN}
Suppose that $T>0$, $\eta\in \mathbb{R}$, $\bu \in \L2 {\HHH}$, the assumptions \ref{ass:1:visco}--\ref{ass:4:uad} are satisfied and the initial data satisfy
$(\bv_0,\ph_0)\in \Vs \times (\Hx5\cap \Hqn)$.
Let the control function $\bu\in\UR $ with $\solst$ being the corresponding state given by Theorem \ref{THM:RES}. For every $\bh \in \L2 \HHH$,
the linearized system \eqref{eq:lin:1}--\eqref{eq:lin:7} admits a unique weak solution $\sollin$ on $[0,T]$ such that
\begin{align*}
	\bw & \in  \L\infty \Hs \cap \L2 \Vs\cap \H1 \Vsp,
	\\
	q &\in  \L2 {H_0},
	\\
	\psi &\in \L\infty {\Hdn} \cap \L2 {\Hx5\cap \Hqn}\cap \H1 H,
	\\
	\theta &\in \L2 V,
	\\
	\w & \in  \L\infty {H} \cap \L2 {\Hx3\cap \Hdn}.
\end{align*}
Moreover, we have
\begin{align*}
& \<\dt \bw,\boldsymbol{z}>_\Vs
	+ \int_\Omega (\bw \cdot \nabla ) \bv \cdot \boldsymbol{z}\,\mathrm{d}x
	+ \int_\Omega  (\bv \cdot \nabla ) \bw \cdot \boldsymbol{z}\,\mathrm{d}x
	+ 2 \int_\Omega   \nu' (\ph) \psi D \bv : \nabla \boldsymbol{z}\,\mathrm{d}x
\\
&\qquad
	 + 2 \int_\Omega \nu (\ph) D \bw  : \nabla \boldsymbol{z}\,\mathrm{d}x
\\
&\quad
= \int_\Omega \theta \nabla \ph \cdot \boldsymbol{z}\,\mathrm{d}x
	+ \int_\Omega  \mu \nabla \psi \cdot \boldsymbol{z}\,\mathrm{d}x
	+ \int_\Omega \bh \cdot \boldsymbol{z}\,\mathrm{d}x,
\end{align*}
and
\begin{align*}
& \int_\Omega \partial_t\psi \xi \,\mathrm{d}x
	+ \int_\Omega (\bw \cdot \nabla \ph)\xi \,\mathrm{d}x
	+ \int_\Omega (\bv \cdot \nabla \psi) \xi \,\mathrm{d}x
    + \int_\Omega m'(\ph)\psi\nabla \mu \cdot \nabla \xi \,\mathrm{d}x
    \\
    &\quad + \int_\Omega m(\ph)\nabla \theta\cdot \nabla \xi \,\mathrm{d}x
	=0,
\end{align*}
almost everywhere in $(0,T)$ for all test functions $\boldsymbol{z}\in \Vs$ and $\xi \in V$, the equations \eqref{eq:lin:4}--\eqref{eq:lin:5} are satisfied almost everywhere in $Q$, the boundary conditions in \eqref{eq:lin:6} except $\dn \theta =0$ are satisfied almost everywhere on $\Sigma$, and $\bw(0)=\0$ in $\Hs $, $\psi(0)=0$ in $H$. Here, a similar consideration can be made as in Remark \ref{RMK:pressure} about the pressure $q$.
\end{proposition}
\begin{proof}[Proof of Proposition \ref{THM:LIN}]
The subsequent estimates are presented in a formal way, while they find direct justifications through a suitable Faedo--Galerkin approximation scheme. Given the conventional nature of this technique, we choose to proceed formally and avoid unnecessary technicalities.

First, testing \eqref{eq:lin:3} by $1$, using integration by parts over $\Omega$ and then integration on $[0,t]\subset [0,T]$, we get
\begin{align}
\int_\Omega \psi(t)\,\mathrm{d}x = \int_\Omega \psi(0)\,\mathrm{d}x =0,\quad \forall\, t\in [0,T].
\label{mean}
\end{align}
Next, inspired by the argument in \cite{WY}, we express the chemical potential $\theta$ as a function of the other variables, thus facilitating the progression of calculations. More precisely, through a combined utilization of equations \eqref{eq:lin:4} and \eqref{eq:lin:5}, we observe that
\begin{align}
	\theta & =
	\label{lin:id}
	\Delta^2 \psi + H(\psi; \ph,\omega),
\end{align}
where
\begin{align}
	\non
	H(\psi; \ph,\omega)&
	=  - \Delta (  f'(\ph) \psi) + f''(\ph) \psi \omega
	+ (f'(\ph) + \eta )(- \Delta \psi + f'(\ph) \psi)
	\\ &   \non
	= - f^{(3)}(\ph)|\nabla \ph|^2 \psi - f''(\ph)\Delta \ph \psi - 2 f''(\ph)\nabla  \ph \cdot \nabla \psi
	- 2f'(\ph)\Delta\psi
	\\ & \quad 	\label{def:H:lin}
	+ f''(\ph)\psi \omega
	+ |f'(\ph)|^2\psi
	- \eta \Delta \psi
	+ \eta  f'(\ph)\psi.
\end{align}
Let us assume, for the moment, that the following estimate holds:
\begin{align}
	\label{control:H}
	\norma{H(\psi; \ph,\omega)}_{V}
	\leq C \norma{\nabla \Delta\psi}
	+C \norma{\psi}_{\Hx2},
\end{align}
where $C$ is a positive constant depending on $K_1$ (see \eqref{control:reg}).
Then, we proceed by testing \eqref{eq:lin:1} with $\bw$, \eqref{eq:lin:3} with $\Delta^2 \psi$, adding the resultants together and using the identity \eqref{lin:id} to obtain
\begin{align*}
	& \frac 12 \frac{\mathrm{d}}{\mathrm{d}t} \big(\norma{\bw}^2 + \norma{\Delta\psi}^2\big)
	+ 2\iO \nu (\ph) D \bw: \nabla \bw\,\mathrm{d}x	
	+ \iO m(\ph) |\nabla \Delta^2\psi|^2 \,\mathrm{d}x
	\\ & \quad
	=
	-\iO (\bw \cdot \nabla ) \bv\cdot \bw \,\mathrm{d}x	
	- 2\iO \nu' (\ph) \psi D \bv : \nabla \bw \,\mathrm{d}x	
    + \iO \theta \nabla \ph  \cdot \bw \,\mathrm{d}x
	\\ & \qquad
	+ \iO \mu \nabla \psi  \cdot \bw \,\mathrm{d}x
	+ \iO \bh \cdot \bw \,\mathrm{d}x
	- \iO m(\ph) \nabla H(\psi; \ph,\omega) \cdot \nabla \Delta^2\psi\,\mathrm{d}x
	\\ & \qquad
	- \iO  \bw \cdot \nabla \ph\Delta^2\psi \,\mathrm{d}x
	- \iO  \bv \cdot \nabla \psi \Delta^2\psi \,\mathrm{d}x
	- \iO m'(\ph)\psi\nabla \mu \cdot \nabla \Delta^2\psi \,\mathrm{d}x
	\\ & \quad
	=:\sum_{i=1}^{9}I_i.
\end{align*}
Using \eqref{app-e8}, \eqref{app-e11}, \eqref{control:H}, the Sobolev embedding $V \emb \Lx \gamma$ in two dimensions, for any $\gamma>1$, along with the \Holder, \Lady, Poincar\'{e}, Poincar\'{e}--Wirtinger, Korn and Young inequalities, we obtain, for any $\delta > 0$, that
\begin{align*}
	I_1 & \leq
	\norma{\bw}_4\norma{\nabla \bv}_4\norma{\bw}
	\\ & \leq
	 C \norma{\bw}^2
    + C \norma{\bw} \norma{\nabla \bw} \norma{\bv}_{\Vs}\norma{ \bv}_{\Ws}  	
    \\ & \leq
	 C \norma{\bw}^2
	+ \delta \norma{D \bw}^2
    + \cd\norma{\bv}_{\Vs}^2\norma{ \bv}_{\Ws}^2\norma{\bw}^2
	\\ & \leq
	 \delta  \norma{D \bw}^2
	 + \cd \norma{\bv}_{\Ws}^2\norma{\bw}^2
     + C \norma{\bw}^2,
\end{align*}
\begin{align*}
	 I_2 & \leq
	 2 \norma{\nu'(\ph)}_\infty \norma{\psi}_\infty\norma{ D \bv}\norma{ D \bw}  \leq
	 \delta \norma{D \bw}^2  + \cd \norma{\psi}^2_{\Hx2},
\end{align*}
\begin{align*}
 I_3 &	 =
	 \iO (\Delta^2 \psi + H(\psi; \ph,\omega)) \nabla \ph \cdot \bw \,\mathrm{d}x
	\\ &  \leq
	 \norma{\Delta^2 \psi}\norma{ \nabla \ph}_4\norma{\bw}_4
	 + \norma{H(\psi; \ph,\omega)}\norma{ \nabla \ph }_4\norma{\bw}_4
	\\ &  \leq
	 C\norma{\Delta^2 \psi}^2
	 + C \norma{H(\psi; \ph,\omega)}^2
	 + C \norma{\ph}_{\Hx2}^2\norma{\bw}\norma{D\bw}
	\\ &  \leq
	\delta \norma{D \bw}^2
	+ C\norma{\Delta^2 \psi}^2
	+ C \norma{H(\psi; \ph,\omega)}^2
	+ \cd \norma{  \ph}_{\Hx2}^4\norma{\bw}^2
	\\ & \leq
	\delta \norma{D \bw}^2
	+ C \norma{\Delta^2 \psi}^2+ C \norma{\nabla \Delta  \psi}^2	+C \norma{\psi}^2_{\Hx2}
	+ \cd \norma{\bw}^2
    \\
    & \leq
	\delta \norma{D \bw}^2
	+ \delta \norma{\nabla \Delta^2 \psi}^2	+ \cd \norma{\psi}^2_{\Hx2} +	\cd \norma{\bw}^2,
\end{align*}
\begin{align*}
	 I_4 &
	 \leq
	 \norma{\mu}_4\norma{\nabla \psi}_4\norma{\bw}
	 \leq
	 \norma{\mu}_V^2\norma{\bw}^2
	 + C \norma{\psi}_{\Hx2}^2,
\end{align*}
\begin{align*}
	 I_5 & \leq
	 \norma{\bh}\norma{\bw}
	 \leq
	  \frac12  \norma{\bh}^2+  \frac 12  \norma{\bw}^2,
\end{align*}
\begin{align*}
	 I_6 & \leq
	 \norma{m(\ph)}_\infty \norma{\nabla H(\psi; \ph,\omega)}\norma{\nabla \Delta ^2\psi}
	 \\
     & \leq \delta\norma{\nabla \Delta ^2\psi}^2
	 + \cd { \norma{\nabla H(\psi; \ph,\omega)}}^2
	 \\
     & \leq \delta\norma{\nabla \Delta ^2\psi}^2 +
             \cd \big(  \norma{\nabla \Delta\psi}^2
	 + \norma{\psi}_{\Hx2}^2\big)
     \\
	 & \leq 2 \delta\norma{\nabla \Delta ^2\psi}^2
	 + \cd \norma{  \psi}^2_{\Hx2},
\end{align*}
\begin{align*}
	 I_7+I_8 & \leq
	 \big(\norma{\bw}_4\norma{\nabla \ph}_4 +\norma{\bv}_4\norma{\nabla \psi}_4\big)\norma{\Delta ^2\psi}	
     \\&
     \leq C\big(\norma{D\bw}^\frac12\norma{\bw}^\frac12\norma{\ph}_{H^2(\Omega)} +\norma{D\bv} \norma{\psi}_{H^2(\Omega)}\big) \norma{\nabla \Delta ^2\psi}^\frac23\norma{\Delta\psi}^\frac13
	 \\ &
	 \leq \delta\norma{\nabla \Delta ^2\psi}^2
	 +\delta\norma{D\bw}^2
	 + \cd \norma{\bw}^2
	 + \cd \norma{\psi}_{\Hx2}^2,
\end{align*}
\begin{align*}
	 I_{9} & \leq
    \norma{m'(\ph)}_\infty\norma{\psi}_4\norma{\nabla \mu}_4
	  \norma{\nabla \Delta ^2\psi}
	 \leq
     \delta\norma{\nabla \Delta ^2\psi}^2
	  + \cd \norma{\mu}^2_{\Hx2}\norma{\psi}_V^2.
\end{align*}
From \eqref{mean} and the Poincar\'{e}--Wirtinger inequality, we find
$$
\|\psi\|\leq C\|\nabla \psi\|\leq C\|\Delta\psi\|^{\frac12}\|\psi\|^{\frac12},
$$
which combined with the elliptic estimate and interpolation yields
\begin{align*}
\|\psi\|_{H^2(\Omega)} &\leq C(\|\Delta\psi\|+\|\psi\|) \leq C\|\Delta \psi\|,\\
\|\psi\|_{H^5(\Omega)} &\leq C(\|\Delta^2\psi\|_V+\|\psi\|) \leq C\big(\|\nabla \Delta^2 \psi\| + \|\Delta \psi\|\big).
\end{align*}
Then recalling \eqref{control:reg}, \ref{ass:1:visco}, \ref{ass:2:mob}, and adjusting $\delta\in (0,1)$ sufficiently small, we infer that
\begin{align*}
	& \frac 12 \frac{\mathrm{d}}{\mathrm{d}t} \big(\norma{\bw}^2
    + \norma{\Delta\psi}^2\big)
	+ \frac {\nu_*}2 \iO  |D \bw |^2\,\mathrm{d}x	
	+ \frac {m_*}2\iO |\nabla \Delta^2\psi|^2\,\mathrm{d}x \\
    & \quad
	\leq C \big(1 + \norma{\bv}^2_{\Ws} + \norma{\mu}_{\Hx2}^2 \big)(\norma{\bw}^2 + \norma{\Delta\psi}^2 ) +  \frac 12 \norma{\bh}^2.
\end{align*}
Thanks to Theorem \ref{THM:RES}, it holds
\begin{align*}
	t \mapsto \big(1 + \norma{\bv(t)}^2_{\Ws} + \norma{\mu(t)}_{\Hx2}^2 \big) \in L^1(0,T).
\end{align*}
Thus, integrating over time and using the initial conditions in \eqref{eq:lin:7}, we infer from Gronwall's lemma that
\begin{align}
	\norma{\bw}_{\L\infty \Hs \cap \L2 \Vs}
	+ \norma{\psi}_{\L\infty {\Hx2} \cap \L2 {\Hx5}}
	\leq C\|\bh\|_{\L2{\HHH}}.
\label{wp-es1}
\end{align}
Furthermore, by comparison in equations \eqref{eq:lin:1}, \eqref{eq:lin:3},  \eqref{eq:lin:4}, and \eqref{eq:lin:5}, it is a standard matter to derive the remaining bounds
\begin{align}
	& \norma{\dt \bw}_{\L2 \Vsp}
	+ \norma{\dt \psi}_{\L2 \Vp}
	+ \norma{q}_{\L2 H} + \norma{\theta}_{\L2 V}\notag
	\\ & \quad
	+\norma{\w}_{\L\infty H \cap \L2 {\Hx3}}\leq C\|\bh\|_{\L2{\HHH}}.
\non
\end{align}

We are left to verify the claimed estimate \eqref{control:H}.
To this end, we deduce from the expression \eqref{def:H:lin} and the estimate \eqref{control:reg} that
\begin{align*}
	 \norma{H(\psi; \ph,\omega)}
	& \leq \|f^{(3)}(\ph)\|_\infty\norma{\nabla \ph}^2_\infty \norma{\psi}
     + \|f''(\ph)\|_\infty \norma{\Delta \ph}_4 \norma{\psi}_4
	+  2\|f''(\ph)\|_\infty  \norma{\nabla  \ph}_\infty\norma{\nabla \psi}
	\\ & \quad +  2\|f'(\ph)\|_\infty \norma{\Delta\psi}
	+ \|f''(\ph)\|_\infty \norma{\psi}_4\norma{\omega}_4
	+ \|f'(\ph)\|_\infty^2 \norma{\psi}
	\\
    &\quad  + |\eta|\norma{ \Delta \psi }
	+ |\eta|\|f'(\ph)\|_\infty \norma{\psi} \\
    & \leq
	C \norma{\psi}_{\Hx2}.
\end{align*}
A direct calculation yields
\begin{align}
	\non
	\nabla H(\psi; \ph,\omega)& =
	- f^{(4)}(\ph)\nabla \ph |\nabla \ph|^2  \, \psi
	- 2 f^{(3)}(\ph)\nabla^2 \ph \nabla \ph  \psi
    - f^{(3)}(\ph)|\nabla \ph|^2 \nabla \psi
	\\ & \quad  \non
	- f^{(3)}(\ph)\nabla \ph \Delta \ph \psi
    - f''(\ph)\nabla \Delta \ph \psi
	- f''(\ph)\Delta \ph \nabla \psi
	\\ & \quad \non
	- 2 f^{(3)}(\ph) \nabla  \ph( \nabla  \ph\cdot \nabla \psi)
	- 2 f''(\ph)\nabla^2 \ph  \nabla \psi
	- 2 f''(\ph) \nabla^2\psi \nabla  \ph
	\\ & \quad \non
	- 2f''(\ph)\nabla \ph\Delta\psi
	- 2f'(\ph)\nabla \Delta\psi
    + f^{(3)}(\ph)\nabla \ph\psi \omega
	+ f''(\ph)\nabla \psi \omega
	\\ & \quad \non	
	+ f''(\ph)\psi \nabla \omega
	+ 2f'(\ph)  f''(\ph)\nabla \ph\psi
    + |f'(\ph)|^2\nabla \psi
	\\ & \quad 	
	- \eta \nabla \Delta \psi
	+ \eta  f''(\ph)\nabla \ph\psi
	+ \eta  f'(\ph)\nabla \psi,
\label{def:nablaH:lin}
\end{align}
where the notation $\nabla^2 \ph$ (or $\nabla^2 \psi$) for the Hessian matrix of $\ph$ (or $\psi$) have been employed.
Using \eqref{def:nablaH:lin}, along with the estimate \eqref{control:reg} and the Sobolev embedding theorem, we can deduce that
\begin{align*}
	 \norma{\nabla H(\psi; \ph,\omega)}
	&\leq
	C \norma{\nabla \ph}^3_\infty \norma{\psi}
	+C  \norma{ \ph}_{W^{2,4}(\Omega)} \norma{\nabla \ph}_\infty  \norma{\psi}_4
	+	C \norma{\nabla \ph}^2_\infty\|\nabla \psi\|
	\\ & \quad
	+C \norma{\nabla \ph}_\infty \norma{ \Delta \ph}_4 \norma{\psi}_4
	+ C \norma{\nabla \Delta \ph}_4 \norma{\psi }_4
	+ C\norma{\Delta \ph}_4\norma{ \nabla \psi}_4
	\\ & \quad
	+C \norma{\nabla  \ph}_\infty^2 \norma{\nabla \psi}
	+ C\norma{\ph}_{W^{2,4}(\Omega)}\norma{\nabla \psi }_4
	+ C\norma{\psi}_{H^2(\Omega)}\norma{\nabla  \ph}_\infty
	\\ & \quad
	+ C \norma{\nabla \ph}_\infty \norma{\Delta\psi}
	+C \norma{\nabla \Delta\psi}
	+ C\norma{\nabla \ph}_\infty\norma{\psi}_4\norma{ \omega}_4
	\\ & \quad
	+ C\norma{\nabla \psi}_4\norma{ \omega}_4
	+ C \norma{\psi }_4\norma{\nabla \omega}_4
    + C\norma{\nabla \ph }_\infty \norma {\psi}
	+ C\norma{\nabla \psi}
	\\ & \quad
	+ C\norma{\nabla \Delta \psi }
	+ C\norma{\nabla \ph}_\infty\norma{\psi}
	+ C\norma{\nabla \psi}
	\\ &
	\leq
	C \norma{\nabla \Delta \psi}
	+
	C \norma{\psi}_{\Hx2}.
\end{align*}
Combining the estimates on $H(\psi; \ph,\omega)$ and its gradient, we can conclude the claim \eqref{control:H}.

The \emph{a priori} estimates obtained above enable us to prove the existence of a weak solution to problem \eqref{eq:lin:1}--\eqref{eq:lin:7} on $[0,T]$ by using the Faedo--Galerkin method and the compactness argument. Moreover, the uniqueness of solutions follows straightforwardly from the linearity of the system.
\end{proof}

Proposition \ref{THM:LIN} paves the way for proving Theorem \ref{THM:FRE} that establishes the \Fre\ differentiability of the control-to-state operator $\Sol$.

\begin{proof}[Proof of Theorem \ref{THM:FRE}]
To show that the solution operator $\Sol$ is \Fre\ differentiable between $\UR$
and the space $\ZZ$ and that, given $\bu\in\UR$, its \Fre\ derivative $D\Sol(\bu)\in {\cal L}(\UR,\ZZ)$ is characterized by the linear operator
that maps any $\bh \in \UR$ into the third component $\psi$ of the solution $\sollin$ to the linearized problem \eqref{eq:lin:1}--\eqref{eq:lin:7}
associated with $\bu$ and the variation $\bh$, we are going to demonstrate that
\begin{equation}
\frac{\norma{\Sol(\bu+\bh) - \Sol(\bu) - \psi  }_{\ZZ}}{\norma{\bh}_{\UR}} \to 0
\quad \text{as}
\quad  \norma{\bh}_{\UR} \to 0.
 \label{frechet}
\end{equation}
 Without loss of generality, we tacitly assume that the sum $\bu+\bh$ satisfies $\bu+\bh\in \UR$. This can be achieved by requiring that the magnitude of the increment $\bh$ is small enough.
 Let us define
 $$
 \Sol(\bu+\bh)=:\phh\quad \text{and}\quad \Sol(\bu)=:\ph
 $$
 the third components of the solutions to the state system \Sys\ corresponding to the given controls $\bu+\bh$ and $\bu$, which we denote by $\solsth$ and $\solst$, respectively. It is worth noting that the uniform bounds in Theorem \ref{THM:RES} are true for both solutions as well as for the corresponding differences.
 Namely, we see that both solutions fulfill \eqref{control:reg} (with bounds depending on $R$) and the following stability estimate
\begin{align}	
		&  \non
		\norma{\bvh-\bv}_{L^\infty(0,t;\Hs)\cap L^2(0,t;\Vs)}
		+ \norma{\phh-\ph}_{L^\infty(0,t;\Hx2) \cap L^2(0,t;\Hx5)}
		\\ & \qquad
		+ \norma{\muh-\mu}_{L^2(0,t;V)}
		+ \norma{\omh-\omega}_{L^\infty(0,t; H) \cap L^2(0,t;\Hx3)}\notag 	\\
		&\quad \leq K_2\norma{\boldsymbol{P}\bh}_{L^2(0,T;\boldsymbol{V}^*_\sigma)},
\quad \forall\, t\in (0,T].
       \label{cd:fre:est}
\end{align}

In order to prove \eqref{frechet}, we introduce the following auxiliary variables
\begin{align*}
	\bs & := \bvh - \bv - \bw,
     \quad
	r:= p^\bh - p -q,\\
     \quad
	\xi &:= \phh - \ph - \psi,
	\quad
	\tau:= \muh-\mu-\theta,
	\quad
	\zeta := \omh - \omega - \w.	
\end{align*}
Then our aim is to show that
\begin{align}
	\label{fre:thesis}
	\norma{\Sol(\bu+\bh) - \Sol(\bu) - \psi}_{\ZZ}
= \norma{\xi}_\ZZ \leq C \norma{\boldsymbol{P}\bh}^2_{\L2 {\boldsymbol{V}^*_\sigma}},
\end{align}
which readily entails \eqref{frechet}, because $\norma{\boldsymbol{P}\bh}_{\L2 {\boldsymbol{V}^*_\sigma}}
\leq \norma{\boldsymbol{P}\bh}_{\L2 {\boldsymbol{H}}}
\leq \norma{\bh}_{\L2 {\boldsymbol{H}}}$.
To this end, upon comparing \Sys\ with \eqref{eq:lin:1}--\eqref{eq:lin:7}, we realize that $(\bs,r,\xi,\tau,\zeta)$ solves the following system
\begin{alignat}{2}
	\non
	& \dt \bs
	+ (\bv \cdot \nabla ) \bs
	+((\bvh -\bv) \cdot \nabla )(\bvh-\bv)
	+ (\bs \cdot \nabla ) \bv
	- \div (2 \nu (\ph) D \bs)
    + \nabla r
	&&
	\\ & \quad 	\label{eq:fre:1}	
	=
	\div {\boldsymbol \Lambda}_1	+  \mu \nabla \xi
	 + (\muh - \mu) \nabla (\phh-\ph)
	 +\tau \nabla \ph
	\qquad && \text{in $Q$},
	\\
	\label{eq:fre:2}
	& \div \bs =0
	\qquad && \text{in $Q$},
	\\
	\label{eq:fre:3}
	& \dt \xi
	+  \bv \cdot \nabla \xi
	+ \bs \cdot \nabla \ph
	- \div (m(\ph)\nabla \tau)
	+ \Lambda_2+ \div {\boldsymbol \Lambda}_3	=0
	\qquad && \text{in $Q$},
	\\
	\label{eq:fre:4}
	& \tau =
	- \Delta \zeta
	+ f'(\ph)\zeta
	+ \Lambda_4
	+ \eta \zeta
	\qquad && \text{in $Q$},
	\\
	\label{eq:fre:5}
	& \zeta =
	-\Delta \xi
	+ f(\phh)- f(\ph) - f'(\ph)\psi
	\qquad && \text{in $Q$},
		\\
	\label{eq:fre:6}
	& \bs = {\boldsymbol 0},
	\quad
	\dn \xi = \dn \Delta \xi = \dn \tau =0
	\qquad && \text{on $\Sigma$},
	\\
	\label{eq:fre:7}
	& \bs(0)=\0,
	\quad
	\xi(0)=0
	\qquad && \text{in $\Omega$},
\end{alignat}
where the terms ${\boldsymbol \Lambda}_1,\Lambda_2,{\boldsymbol \Lambda}_3$, and $\Lambda_4$ are defined by
\begin{align}
	\label{def:L1}
	{\boldsymbol \Lambda}_1 & :=   2 (\nu (\phh)-\nu (\ph))(D \bvh-D \bv)
	+ 2 (\nu(\phh)-\nu(\ph) -\nu' (\ph) \psi) D \bv ,
	\\
	\label{def:L2}
	\Lambda_2 & :=	(\bvh-\bv)\cdot\nabla(\phh-\ph),
	\\
	\label{def:L3}
	{\boldsymbol \Lambda}_3 & :=  (m(\phh)-m(\ph))\nabla (\muh-\mu)
	+  (m(\phh)-m(\ph)- m'(\ph)\psi)\nabla \mu,
	\\
	\label{def:L4}
	\Lambda_4 & := (f'(\phh)-f'(\ph))(\omh - \omega)
	+ (f'(\phh)- f'(\ph) - f''(\ph)\psi)\omega.
\end{align}

Before proceeding, let us recall Taylor's formula with integral remainder for a $C^2$ function $g :\erre \to \erre$, that is,
\begin{align*}
	g(\phh)-g(\ph)-g'(\ph)\psi = g'(\ph)\xi + {\cal R}_g (\phh-\ph)^2,
\end{align*}
with
$$
{\cal R}_g = \int_0^1 \int_0^1 g''\big(sz \phh+(1-sz)\ph)\big)z \, \mathrm{d}s\mathrm{d}z.
$$
The remainder term ${\cal R}_g$ is uniformly bounded due to the uniform boundedness of $\phh$, $\ph$ and because the integration variables $s$, $z$ belong to $[0,1]$. Thus, there exists a positive constant $C$ such that $\norma{{\cal R}_g}_{L^\infty(Q)} \leq C$.
In the subsequent analysis, we will employ this formula with $g$ being any of $g\in\{\nu,m,f,f'\}$, cf. \eqref{eq:fre:5}, \eqref{def:L1}, \eqref{def:L3}, and \eqref{def:L4}.

From equation \eqref{eq:fre:3}, it easily follows that
\begin{align}
\int_\Omega \xi(t)\,\mathrm{d}x=\int_\Omega \xi(0)\,\mathrm{d}x=0,\quad
\forall\, t\in [0,T].
\label{mean-xi}
\end{align}
Next, in view of \eqref{eq:fre:4} and \eqref{eq:fre:5}, we infer that
\begin{align}
	\non
	\tau & = \Delta^2 \xi
	- \Delta  (f(\phh)- f(\ph) - f'(\ph)\psi)
	- f'(\ph)\Delta \xi
	\\ & \quad \non
	+ f'(\ph) (f(\phh)- f(\ph) - f'(\ph)\psi)
	+ \Lambda_4 + \eta
	(-\Delta \xi
	+ f(\phh)- f(\ph) - f'(\ph)\psi)
	\\ & =: \label{eq:fre:combined}
	\Delta^2 \xi
	+ G(\xi;\phh,\ph,\psi).
\end{align}
Let us claim that the following estimate holds:
\begin{align}
&	\norma{G(\xi;\phh,\ph,\psi)}_{L^2(0,t;V)} \non\\
&\quad  \leq C\norma{\nabla \Delta \xi}_{L^2(0,t;\HHH)} + C \norma{\xi}_{L^2(0,t;\Hx2)}
	+ C \norma{\boldsymbol{\bh}}_{L^2(0,T;\boldsymbol{V}_\sigma^*)}^2,\quad \forall\, t\in [0,T].
	\label{eq:fre:G:H1}
\end{align}
Then, we test \eqref{eq:fre:1} by $\bs$, \eqref{eq:fre:3} by $\Delta^2 \xi$, make use of the reformulation \eqref{eq:fre:combined} and add the resulting identities together. Upon rearranging terms and integrating by parts, we obtain
\begin{align*}
	& \frac 12 \frac{\mathrm{d}}{\mathrm{d}t} \big(\norma{\bs}^2
     +   \norma{\Delta\xi}^2\big)
	+  2 \iO \nu (\ph)  D \bs: \nabla \bs  \,\mathrm{d}x
	+ \iO m(\ph) |\nabla \Delta^2\xi|^2 \,\mathrm{d}x
	\\ & \quad
	=
	 - \iO  \boldsymbol \Lambda_1	 : \nabla \bs \,\mathrm{d}x
	+\iO  \mu \nabla \xi \cdot \bs\,\mathrm{d}x
	+ \iO (\muh - \mu) \nabla (\phh-\ph) \cdot \bs\,\mathrm{d}x
	\\ & \qquad
	+ \iO \tau \nabla \ph \cdot \bs\,\mathrm{d}x
	- \iO ((\bvh -\bv) \cdot \nabla )(\bvh-\bv) \cdot \bs\,\mathrm{d}x
		\\ & \qquad
	- \iO  (\bs \cdot \nabla ) \bv \cdot \bs\,\mathrm{d}x
	-\iO \bv \cdot \nabla \xi \Delta^2\xi\,\mathrm{d}x
	- \iO \bs \cdot \nabla \ph \Delta^2\xi\,\mathrm{d}x
		\\ & \qquad
	- \iO m(\ph) \nabla G(\xi;\phh,\ph,\psi) \cdot \nabla \Delta^2\xi\,\mathrm{d}x
	- \iO \Lambda_2\Delta^2\xi\,\mathrm{d}x
	+ \iO {\boldsymbol \Lambda}_3 \cdot \nabla \Delta^2\xi\,\mathrm{d}x\\
    &\quad
	=: \sum_{i=1}^{11} I_i.
\end{align*}
Using the interpolation inequality \eqref{inter}, the estimate \eqref{control:reg}, Taylor's formula mentioned above, the Young, Korn, and \Holder\ inequalities, we obtain, for any $\delta >0$, that
\begin{align}
	I_1 &
	=
	- 2 \iO (\nu (\phh)-\nu (\ph))(D \bvh-D \bv) : \nabla  \bs\,\mathrm{d}x
	- 2 \iO (\nu(\phh)-\nu(\ph) -\nu' (\ph) \psi) D \bv : \nabla \bs\,\mathrm{d}x
    \notag
	\\ & \leq
	C \norma{\phh-\ph}_\infty\norma{D \bvh-D \bv}\norma{D \bs} +
	C (\norma{\xi}_\infty+\norma{\phh-\ph}^2_\infty) \norma{D \bv}
	\norma{D \bs} \notag
	\\ & \leq
	\delta \norma{D\bs}^2
	+ \cd \norma{\phh-\ph}^2_{H^2(\Omega)}\norma{\bvh-\bv}^2_{\Vs}
	+ \cd \norma{\xi}_{H^2(\Omega)}^2
	+ \cd  \norma{\phh-\ph}^4_{H^2(\Omega)},
    \notag
\end{align}
\begin{align}
	I_2 & \leq
	\norma{\mu}_{4}\norma{\nabla \xi }_{4}\norma{\bs}
	\leq  C\norma{\bs}^2 +
	C \norma{\mu}^2_{V}\norma{\xi }^2_{\Hx2},
    \notag
\end{align}
\begin{align*}
	I_3 & \leq
	\norma{\muh - \mu}_4 \norma{\nabla (\phh-\ph)}_{4} \norma{\bs}
        \leq
	C \norma{\bs}^2
	+ C \norma{\muh - \mu}_{V}^2 \norma{\phh-\ph}_{\Hx2}^2 , 
\end{align*}
\begin{align*}
	I_4 &
	= \iO (\Delta^2 \xi
	+ G(\xi;\phh,\ph,\psi)) \nabla \ph \cdot \bs\,\mathrm{d}x \\
    & 	\leq C \big(\norma{\Delta^2 \xi } + \norma{G(\xi;\phh,\ph,\psi) }\big)\norma{\nabla \ph }_4\norma{\bs}_4 \\
    & \leq
	\delta \norma{\nabla \bs}^2
	+ \cd \norma{\Delta^2 \xi }^2\norma{\ph }^2_{\Hx2}
    + \cd \norma{G(\xi;\phh,\ph,\psi)}^2\norma{\ph }^2_{\Hx2} \\
    &\leq  2\delta \norma{D\bs}^2 + \delta \norma{\nabla \Delta^2\xi}^2
    + \cd \norma{ \ph}_{\Hx2}^6 \norma{ \xi}_{\Hx2}^2
    + \cd \norma{G(\xi;\phh,\ph,\psi)}^2\norma{\ph }^2_{\Hx2},
\end{align*}
\begin{align*}
	{I_5}
    &= \iO ((\bvh -\bv) \cdot \nabla ) \bs \cdot (\bvh-\bv) \,\mathrm{d}x\\
    & \leq
	\norma{\bvh -\bv}_4^2\norma{\nabla \bs }\\
	& \leq \delta \norma{\nabla \bs}^2
	+\cd \norma{\bvh -\bv}_4^4
	\\ & \leq
	2\delta \norma{D\bs}^2
	+\cd \norma{\bvh -\bv}^2\norma{\bvh-\bv}_{\Vs}^{2}, 
\end{align*}
\begin{align*}
	I_6 & \leq
	\norma{\bs}_4\norma{\nabla \bv}_4\norma{\bs}
	\leq \delta \norma{D\bs}^2
	+ \cd \norma{\bv }_{\Ws}^2 \norma{\bs}^2 ,
\end{align*}
\begin{align*}
	I_7 & \leq
	\norma{\bv }_4\norma{\nabla \xi }_4\norma{\Delta^2\xi}
	\leq C\norma{\bv }_{\Vs}^2\norma{\xi }_{\Hx2}^2
	+C \norma{\Delta^2\xi}^2\\
    &\leq \delta \norma{\nabla \Delta^2\xi}^2 + C_\delta\norma{\xi }_{\Hx2}^2 + C\norma{\bv }_{\Vs}^2\norma{\xi }_{\Hx2}^2,
\end{align*}
\begin{align*}
	I_8 & \leq
	\norma{ \bs }_4\norma{ \nabla \ph}_4\norma{ \Delta^2\xi}
	\leq \delta \norma{D\bs}^2
	+ \cd \norma{ \ph}_{\Hx2}^2\norma{ \Delta^2\xi}^2
   \\
   &\leq \delta \norma{D\bs}^2  + \delta \norma{\nabla \Delta^2\xi}^2
   + \cd \norma{ \ph}_{\Hx2}^6\norma{ \xi}_{\Hx2}^2,
\end{align*}
\begin{align*}
	I_{9} & \leq
	\|m(\ph)\|_{L^\infty(\Omega)} \norma{\nabla G(\xi;\phh,\ph,\psi)}\norma{\nabla \Delta^2\xi}
    \\
    & \leq \delta \norma{\nabla \Delta^2\xi}^2
         + \cd  \norma{\nabla G(\xi;\phh,\ph,\psi)}^2,
\end{align*}
\begin{align*}
	I_{10} & \leq
	\norma{\bvh-\bv}_4\norma{\nabla \diff}_4
	\norma{ \Delta^2\xi} \\
    &\leq
	C	\norma{\bvh-\bv}_{\Vs}^2\norma{ \phh-\ph}_{\Hx2}^2
	+ C \norma{ \Delta^2\xi}^2\\
    &\leq \delta \norma{\nabla \Delta^2\xi}^2 + C	\norma{\bvh-\bv}_{\Vs}^2\norma{ \phh-\ph}_{\Hx2}^2
    + \cd \norma{ \xi}_{\Hx2}^2, 
\end{align*}
\begin{align*}
	I_{11} & =
	\iO (m(\phh)-m(\ph))\nabla (\muh-\mu)
	 \cdot \nabla \Delta^2\xi\,\mathrm{d}x \\
   &\quad
	 + \iO\big( m'(\ph)\xi + \rem_m \diff^2\big)\nabla \mu \cdot \nabla \Delta^2\xi \,\mathrm{d}x
	 \\ &
	 \leq  \delta \norma{\nabla \Delta^2\xi}^2
	 + \cd \norma{\phh-\ph}_{\infty}^2 \norma{\muh-\mu}_{V}^2
	 + \cd \norma{\xi}_\infty^2\norma{\mu}_V^2
	 + \cd \norma{\phh-\ph}_{\infty}^4\norma{\mu}_V^2\\
&\leq \delta \norma{\nabla \Delta^2\xi}^2
	 + \cd \norma{\phh-\ph}_{H^2(\Omega)}^2 \norma{\muh-\mu}_{V}^2
	 + \cd \norma{\xi}_{H^2(\Omega)}^2\norma{\mu}_V^2
	 \\ & \quad
	 + \cd \norma{\phh-\ph}_{H^2(\Omega)}^4\norma{\mu}_V^2.
\end{align*}
Collecting the above estimates, using assumptions \ref{ass:1:visco}, \ref{ass:2:mob}, the estimate \eqref{control:reg}, the mass conservation property \eqref{mean-xi}, and taking $\delta$ small enough, we obtain
\begin{align*}
	& \frac 12 \frac{\mathrm{d}}{\mathrm{d}t} \big(\norma{\bs}^2
     +   \norma{\Delta\xi}^2\big)
	+ \frac{\nu_*}{2}\iO   |D \bs|^2 \,\mathrm{d}x
	+ \frac{m_*}{2}\iO   |\nabla \Delta^2\xi|^2 \,\mathrm{d}x\\
&\quad \leq C\big(1+\norma{\bv }_{\Ws}^2\big)\big(\norma{\bs}^2+\norma{\Delta\xi}^2\big)
+ C \norma{\bvh -\bv}^2\norma{\bvh-\bv}_{\Vs}^2
\\
&\qquad + C \norma{\phh-\ph}^2_{H^2(\Omega)} \norma{\bvh-\bv}_{\Vs}^2
+ C \norma{\muh - \mu}_{V}^2 \norma{\phh-\ph}^2_{H^2(\Omega)}\\
&\qquad + C \norma{\phh-\ph}^4_{H^2(\Omega)} + C \norma{G(\xi;\phh,\ph,\psi)}_V^2.
\end{align*}
Integrating the above differential inequality over time, from \eqref{control:reg}, \eqref{cd:fre:est}, \eqref{mean-xi}, the claim \eqref{eq:fre:G:H1}, the interpolation inequality, as well as Young's inequality, we can deduce that
\begin{align*}
&\norma{\bs(t)}^2 + \norma{\Delta\xi(t)}^2 +
	  \nu_*\int_0^t  \|D \bs(s)\|^2 \,\mathrm{d}s +  m_*\int_0^t  \|\nabla \Delta^2\xi(s)\|^2 \,\mathrm{d}s\\
&\quad \leq C\int_0^t\big(1+\norma{\bv(s) }_{\Ws}^2\big)\big(\norma{\bs(s)}^2+\norma{\Delta\xi(s)}^2\big)\,\mathrm{d}s\\
& \qquad + C \norma{\boldsymbol{P}\bh}_{L^2(0,t;\boldsymbol{V}_\sigma^*)}^4 +C \int_0^t \|\nabla \Delta\xi(s)\|^2\,\mathrm{d}s\nonumber\\
&\quad \leq  C\int_0^t\big(1+\norma{\bv(s) }_{\Ws}^2\big)\big(\norma{\bs(s)}^2+\norma{\Delta\xi(s)}^2\big)\,\mathrm{d}s\\
&\qquad  + C \norma{\boldsymbol{P}\bh}_{L^2(0,T;\boldsymbol{V}_\sigma^*)}^4 + \frac{m_*}{2} \int_0^t \|\nabla \Delta^2\xi(s)\|^2 \,\mathrm{d}s,\quad \forall\, t\in[0,T].
\end{align*}
Then an application of Gronwall's lemma yields
\begin{align*}
	\|\bs\|_{L^\infty(0,T;\Hs)\cap L^2(0,T;\Vs)}+  \norma{\xi}_{\L\infty{H^2(\Omega)}\cap \L2{H^5(\Omega)}}
	\leq
	C \norma{\boldsymbol{P}\bh}_{L^2(0,T;\boldsymbol{V}_\sigma^*)}^2.
\end{align*}
From this estimate and \eqref{eq:fre:3}, it is standard to deduce that
\begin{align*}
	\norma{\dt \xi}_{\L2 \Vp} \leq C \norma{\boldsymbol{P}\bh}_{L^2(0,T;\boldsymbol{V}_\sigma^*)}^2.
\end{align*}
Collecting the above estimates, we arrive at the conclusion \eqref{fre:thesis}.

It remains to verify the claimed estimate \eqref{eq:fre:G:H1}. We observe that in the expression of $G(\xi;\phh,\ph,\psi)$, it occurs the following critical term
$$
\Delta  (f(\phh)- f(\ph) - f'(\ph)\psi) = \Delta (f'(\ph)\xi + \widehat{\rem}),
$$
where
\begin{align*}
	\widehat{\rem} = \rem_{f}\diff^2= \left(\int_0^1\int_0^1 f''(sz \phh+(1-sz)\ph) z\, \mathrm{d}s\mathrm{d}z \right)\diff^2.
\end{align*}
A direct calculation yields
\begin{align*}
	\nabla \widehat{\rem}
	& = 2 \diff \nabla \diff\int_0^1\int_0^1 f''(sz \phh+(1-sz)\ph) z\, \mathrm{d}s\mathrm{d}z
	\\ & \quad
	+ \diff^2 \int_0^1\int_0^1 f^{(3)}(sz \phh+(1-sz)\ph)(sz \nabla \phh+(1-sz)\nabla \ph)z \,\mathrm{d}s \mathrm{d}z,
\end{align*}
\begin{align*}
	\Delta \widehat{\rem} & =
	2 |\nabla \diff|^2 \int_0^1\int_0^1 f''(sz \phh+(1-sz)\ph) z\, \mathrm{d}s\mathrm{d}z
	\\ & \quad
	+ 2 \diff \Delta \diff \int_0^1\int_0^1 f''(sz \phh+(1-sz)\ph) z\, \mathrm{d}s\mathrm{d}z
	\\ & \quad
	+ 4 \diff \nabla \diff \cdot \Big(\int_0^1\int_0^1 f^{(3)}(sz \phh+(1-sz)\ph)(sz \nabla \phh+(1-sz)\nabla \ph)z \,\mathrm{d}s \mathrm{d}z\Big)
	\\ & \quad
	+ \diff^2\int_0^1\int_0^1 f^{(4)}(sz \phh+(1-sz)\ph)|sz \nabla \phh+(1-sz)\nabla \ph|^2z \,\mathrm{d}s \mathrm{d}z
	\\ & \quad
	+ \diff^2 \int_0^1\int_0^1 f^{(3)}(sz \phh+(1-sz)\ph)(sz \Delta \phh+(1-sz)\Delta \ph)z \,\mathrm{d}s \mathrm{d}z,
\end{align*}
and
\begin{small}
\begin{align*}
	\nabla \Delta \widehat{\rem}
& =
	4  \nabla^2 \diff \nabla \diff \int_0^1\int_0^1 f''(sz \phh+(1-sz)\ph) z\, \mathrm{d}s\mathrm{d}z
	\\ & \quad
	+2 |\nabla \diff |^2\int_0^1\int_0^1 f^{(3)}(sz \phh+(1-sz)\ph)(sz \nabla \phh+(1-sz)\nabla \ph)z \,\mathrm{d}s \mathrm{d}z
	\\ & \quad
	 + 2\nabla \diff \Delta \diff \int_0^1\int_0^1 f''(sz \phh+(1-sz)\ph) z\, \mathrm{d}s\mathrm{d}z
	 \\ & \quad
	+2  \diff  \nabla \Delta \diff \int_0^1\int_0^1 f''(sz \phh+(1-sz)\ph) z\, \mathrm{d}s\mathrm{d}z
	\\ & \quad
	+2 \diff  \Delta \diff \int_0^1\int_0^1 f^{(3)}(sz \phh+(1-sz)\ph)(sz \nabla \phh+(1-sz)\nabla \ph)z \,\mathrm{d}s \mathrm{d}z
	 	\\ & \quad
	+ 4
	\nabla \diff \otimes \nabla \diff
	\Big(\int_0^1\int_0^1 f^{(3)}(sz \phh+(1-sz)\ph)(sz \nabla \phh+(1-sz)\nabla \ph)z \,\mathrm{d}s \mathrm{d}z\Big)
	 	\\ & \quad
	 + 4 \diff \nabla^2 \diff \Big(\int_0^1\int_0^1 f^{(3)}(sz \phh+(1-sz)\ph)(sz \nabla \phh+(1-sz)\nabla \ph)z \,\mathrm{d}s \mathrm{d}z\Big)
	\\ & \quad
	+ 4 \diff \Big(\int_0^1\int_0^1  f^{(4)}(sz \phh+(1-sz)\ph)(sz \nabla \phh+(1-sz)\nabla \ph)\\
&\qquad \quad \times \big(\nabla \diff \cdot(sz \nabla \phh+(1-sz)\nabla \ph)\big) z\, \mathrm{d}s\mathrm{d}z\Big)
	\\ & \quad
	+ 4 \diff \Big(\int_0^1\int_0^1 f^{(3)}(sz \phh+(1-sz)\ph)(sz \nabla^2 \phh+(1-sz)\nabla^2 \ph)z \,\mathrm{d}s \mathrm{d}z\Big)  \nabla \diff
	\\ & \quad
	+ 2 \diff \nabla\diff \Big(\int_0^1\int_0^1 f^{(4)}(sz \phh+(1-sz)\ph)|sz \nabla \phh+(1-sz)\nabla \ph|^2z \,\mathrm{d}s \mathrm{d}z\Big)
	\\ & \quad
+ \diff^2 \Big(\int_0^1\int_0^1  f^{(5)}(sz \phh+(1-sz)\ph)|sz \nabla\phh+(1-sz)\nabla\ph|^2 (sz \nabla\phh+(1-sz)\nabla\ph)z \, \mathrm{d}s \mathrm{d}z\Big)
	\\ & \quad
+ 2 \diff^2 \int_0^1\int_0^1  f^{(4)}(sz \phh+(1-sz)\ph) (sz \nabla^2 \phh+(1-sz)\nabla^2 \ph) (sz \nabla \phh+(1-sz)\nabla \ph)z\, \mathrm{d}s\mathrm{d}z
	\\ & \quad
	+ 2 \diff \nabla \diff \int_0^1\int_0^1 f^{(3)}(sz \phh+(1-sz)\ph)(sz \Delta \phh+(1-sz)\Delta \ph)z \,\mathrm{d}s \mathrm{d}z
	\\ & \quad
	+ \diff^2 \int_0^1\int_0^1 f^{(4)}(sz \phh+(1-sz) \ph) (sz \nabla \phh+(1-sz)\nabla \ph)(sz \Delta  \phh+(1-sz)\Delta \ph)z \, \mathrm{d}s\mathrm{d}z
	\\ & \quad
	+ \diff^2 \int_0^1 \int_0^1 f^{(3)}(sz \phh+(1-sz)\ph)(sz \nabla \Delta \phh+(1-sz)\nabla \Delta \ph) z\,\mathrm{d}s \mathrm{d}z.
\end{align*}
\end{small}
Using the above identities and the estimate \eqref{control:reg}, we can deduce that
\begin{align*}
	|\Delta \widehat{\rem} |
	& \leq C |\nabla \diff|^2
	+C| \phh-\ph||\Delta \diff|
	\\ & \quad
	 + C |\phh-\ph||\nabla \diff| (|\nabla \phh|+|\nabla \ph|)
	 \\ & \quad
	 +C |\phh-\ph|^2(|\nabla \phh|^2+|\nabla \ph|^2+ |\Delta \phh|+|\Delta\ph|),
\end{align*}
and
\begin{align*}
	|\nabla \Delta \widehat{\rem} |
	& \leq C
	|\nabla  \diff|\big(|\nabla^2\diff|+|\Delta\diff|\big)\\
	& \quad + C |\nabla \diff|^2	(|\nabla \phh|+|\nabla \ph|)
	+ C |\phh-\ph||\nabla \Delta \diff|
	\\ & \quad
	+ C |\phh-\ph|\big(|\nabla^2\diff|+ |\Delta \diff|\big)\big(|\nabla \phh|+|\nabla \ph|\big)
	\\ & \quad
	+ C |\phh-\ph||\nabla \diff|
	\big(|\nabla \phh|^2+|\nabla \ph|^2 + |\Delta \phh|+|\Delta \ph|+|\nabla^2 \phh|+|\nabla^2 \ph|\big)
	\\ & \quad
	+ C |\phh-\ph|^2
	\big(|\nabla \phh|^3+|\nabla \ph|^3
	+ |\nabla\Delta \phh|+|\nabla\Delta \ph|\big)
		\\ & \quad
	+ C |\phh-\ph|^2  (|\nabla \phh|+|\nabla \ph|)(|\Delta\phh|+|\Delta \ph|+|\nabla^2\phh|+|\nabla^2 \ph|).
\end{align*}
Then using \Holder's inequality and the stability estimate \eqref{cd:fre:est}, recalling in particular that $\phh $ and $\ph$ are both bounded in $\L\infty {\Hx5}$, we find that for any $t\in [0,T]$, it holds
\begin{align*}
	  \norma{\Delta \widehat{\rem}}_{L^2(0,t;H)}
	& \leq C \norma{\phh-\ph}_{L^\infty(0,t;W^{1,4}(\Omega))}^2
	+C \norma{\phh-\ph}_{L^\infty(Q_t)}\norma{\phh-\ph}_{L^2(0,t;\Hx2)}
	\\ & \quad
	+ C \norma{\phh-\ph}_{L^\infty (Q_t)}\norma{\phh-\ph}_{L^2(0,t; V)}
	+ C \norma{\phh-\ph}_{L^\infty(Q_t)}^2.
\end{align*}

Let us first estimate $G(\xi;\phh,\ph,\psi)$. It follows from the above estimates that
\begin{align}
	& \norma{G(\xi;\phh,\ph,\psi)}_{L^2(0,t; H)}
\notag \\
	&\quad \leq
	\norma{\Delta  (f(\phh)- f(\ph) - f'(\ph)\psi)}_{L^2(0,t; H)}
	+ \norma{f'(\ph)\Delta \xi}_{L^2(0,t; H)}
	\notag \\ & \qquad
		+ \norma{ f'(\ph) (f(\phh)- f(\ph) - f'(\ph)\psi)}_{L^2(0,t; H)}
	+ \norma{\Lambda_4 }_{L^2(0,t; H)}
		+ \norma{\eta \zeta}_{L^2(0,t; H)}
		\notag \\
    & \quad
		=:\sum_{i=1}^{5} I_i,\quad \forall\, t\in [0,T].
\label{es-GL2}
\end{align}
Recall \eqref{def:L4} and the identity
\begin{align*}
	\Delta  (f'(\ph)\xi)
	= f^{(3)}(\ph) |\nabla \ph|^2 \xi
	+ f''(\ph) \Delta\ph \, \xi
	+ 2 f''(\ph)  \nabla \ph \cdot \nabla \xi
	+ f'(\ph) \Delta \xi.
\end{align*}
Then by \eqref{control:reg}, \eqref{cd:fre:est}, the terms on the \rhs\ of \eqref{es-GL2} can be bounded as follows
\begin{align*}
	I_1 & \leq	C \norma{\Delta (f'(\ph)\xi)}_{L^2(0,t; H)}
		+ C \norma{\Delta \widehat{\rem}}_{L^2(0,t; H)}
		\\
        & \leq
		C \norma{\xi}_{L^2(0,t;{\Hx2})}
		+C\norma{\boldsymbol{P}\bh}_{L^2(0,T;\boldsymbol{V}_\sigma^*)}^2,
\end{align*}
\begin{align*}
		I_2+I_5
        & \leq
		C \norma{\xi}_{L^2(0,t;{\Hx2})}
        +C\norma{\boldsymbol{P}\bh}_{L^2(0,T;\boldsymbol{V}_\sigma^*)}^2,
\end{align*}
\begin{align*}
		I_3
		& \leq \norma{ |f'(\ph)|^2\xi}_{L^2(0,t; H)} + \norma{ f'(\ph) \widehat{\rem} }_{L^2(0,t; H)}\\
		& \leq C \norma{\xi}_{L^2(0,t; H)}
		+C\norma{\boldsymbol{P}\bh}_{L^2(0,T;\boldsymbol{V}_\sigma^*)}^2,
\end{align*}
\begin{align*}
		I_4
        &\leq \|(f'(\phh)-f'(\ph))(\omh-\omega)\|_{L^2(0,t; H)}
        + \|(f''(\ph)\xi+ \rem_{f'}(\phh-\ph)^2)\omega\|_{L^2(0,t; H)}
        \\
		& \leq
		C \norma{\phh-\ph}_{L^\infty(Q_t)} \norma{\omh-\omega}_{L^2(0,t; H)}
		+ C \norma{\xi}_{L^2(0,t;{L^\infty(\Omega)})} \norma{\omega}_{L^\infty(0,t;H)}
		\\ & \quad
		+ C \norma{\phh-\ph}_{L^\infty(Q_t)}^2 \norma{\omega}_{L^\infty(0,t;H)}.
\end{align*}
Hence, it holds
\begin{align*}
	& \norma{G(\xi;\phh,\ph,\psi)}_{L^2(0,t; H)}
	 \leq
	C \norma{\xi}_{L^2(0,t;\Hx2)}
	+ C \norma{\boldsymbol{P}\bh}_{L^2(0,T;\boldsymbol{V}_\sigma^*)}^2,
\quad \forall\, t\in [0,T].
\end{align*}
Similarly, to control the gradient term $\nabla G(\xi;\phh,\ph,\psi)$, we observe that
\begin{align}
	& \norma{\nabla G(\xi;\phh,\ph,\psi)}_{L^2(0,t; \HHH)}
\notag \\
	&\quad \leq
	\norma{\nabla \Delta  (f(\phh)- f(\ph) - f'(\ph)\psi)}_{L^2(0,t; \HHH)}
	+ \norma{\nabla (f'(\ph)\Delta \xi)}_{L^2(0,t; \HHH)}
	\notag \\ & \qquad
		+ \norma{\nabla( f'(\ph) (f(\phh)- f(\ph) - f'(\ph)\psi))}_{L^2(0,t; \HHH)}
	+ \norma{\nabla \Lambda_4 }_{L^2(0,t; \HHH)}
		+ \norma{\eta \nabla \zeta}_{L^2(0,t; \HHH)}
		\notag \\
    & \quad
		=:\sum_{i=1}^{5} I_i,\quad \forall\, t\in [0,T].
\label{es-GH1}
\end{align}
With \eqref{control:reg}, we find
\begin{align*}
	&\norma{\nabla \Delta \widehat{\rem}}_{L^2(0,t; \HHH)}\\
	&\quad  \leq
	C \norma{\phh-\ph}_{L^\infty(0,t;W^{1,4}(\Omega))} \norma{\phh-\ph}_{L^2(0,t;{W^{2,4}(\Omega)})}
	+C \norma{\phh-\ph}_{L^\infty(0,t; W^{1,4}(\Omega))}^2
	\\ & \qquad
	+C \norma{\phh-\ph}_{L^\infty (Q_t)}\norma{\phh-\ph}_{L^2(0,t;{\Hx3})}
    + C \norma{\phh-\ph}_{L^\infty (Q_t)} \norma{\phh-\ph}_{L^\infty(0,t;H^2(\Omega))} \\
    &\qquad
	+ C\norma{\phh-\ph}_{L^\infty (Q_t)}^2.
\end{align*}
From \eqref{control:reg}, \eqref{cd:fre:est}, we see  that for any $t\in [0,T]$, the terms on the right-hand side of \eqref{es-GH1} can be estimated as follows
\begin{align*}
	I_1 & \leq	C \norma{\nabla \Delta (f'(\ph)\xi)}_{L^2(0,t; \HHH)}
		+ C \norma{\nabla \Delta \widehat{\rem}}_{L^2(0,t; \HHH)}
		\\
        & \leq C\norma{\nabla \Delta \xi}_{L^2(0,t;\HHH)} +
		C \norma{\xi}_{L^2(0,t;{\Hx2})}
		+C\norma{\boldsymbol{P}\bh}_{L^2(0,T;\boldsymbol{V}_\sigma^*)}^2,
\end{align*}
\begin{align*}
		I_2+I_5
        & \leq C\norma{\nabla \Delta \xi}_{L^2(0,t;\HHH)} +
		C \norma{\xi}_{L^2(0,t;{\Hx2})}
        +C\norma{\boldsymbol{P}\bh}_{L^2(0,T;\boldsymbol{V}_\sigma^*)}^2,
\end{align*}
\begin{align*}
		I_3
		& \leq \norma{\nabla( |f'(\ph)|^2\xi)}_{L^2(0,t; \HHH)} + \norma{ \nabla (f'(\ph) \widehat{\rem}) }_{L^2(0,t; \HHH)}\\
		& \leq C \norma{\xi}_{L^2(0,t; V)}
		+ C  \norma{\phh-\ph}_{L^\infty(Q_t)}^2
        + C  \norma{\phh-\ph}_{L^\infty(Q_t)}\norma{\phh-\ph}_{L^\infty(0,t;V)},
\end{align*}
\begin{align*}
		I_4
       &\leq \|\nabla [(f'(\phh)-f'(\ph))(\omh-\omega)]\|_{L^2(0,t; \HHH)}
        + \|\nabla [(f''(\ph)\xi+  \rem_{f'}(\phh-\ph)^2)\omega]\|_{L^2(0,t; \HHH)}
        \\
		& \leq C \norma{\phh-\ph}_{L^\infty(0,t; W^{1,4}(\Omega))}\norma{\omh-\omega}_{L^2(0,t;L^4(\Omega))}
         + C \norma{ {\phh-\ph}}_{L^\infty(Q_t)}\norma{\omh-\omega}_{L^2(0,t;H)}
         \\
      &\quad +C \norma{ {\phh-\ph}}_{L^\infty(Q_t)}\norma{\nabla (\omh-\omega)}_{L^2(0,t; \HHH)}
      +C \norma{\xi}_{L^2(0,t; L^4(\Omega))}\norma{\omega}_{L^\infty(0,t; L^4(\Omega))}
      \\
         &\quad
         +C \norma{\nabla \xi}_{L^2(0,t;\boldsymbol{L}^4(\Omega))} \norma{\omega}_{L^\infty(0,t; L^4(\Omega))}
          + C \norma{ {\phh-\ph}}_{L^\infty(Q_t)}^2
          \norma{\omega}_{L^2(0,t; L^2(\Omega))}
          \\
         &\quad
         + C \norma{ {\phh-\ph}}_{L^\infty(Q_t)}
         \norma{ \nabla(\phh-\ph)}_{L^\infty(0,t;\boldsymbol{L}^4(\Omega))}
          \norma{\omega}_{L^2(0,t; L^4(\Omega))}
          \notag\\
         &\quad
         +C (\norma{\xi}_{L^2(0,t;L^4(\Omega))} +\norma{ {\phh-\ph}}_{L^\infty(Q_t)}^2 ) \norma{\nabla  \omega}_{L^\infty(0,t; \boldsymbol{L}^4(\Omega))}
         \\
         &\leq  C \norma{\xi}_{L^2(0,t;{\Hx2})}
		+C\norma{\boldsymbol{P}\bh}_{L^2(0,T;\boldsymbol{V}_\sigma^*)}^2,
\end{align*}
for any $t\in [0,T]$.
As a consequence, we obtain
\begin{align*}
	\norma{\nabla &G(\xi;\phh,\ph,\psi)}_{L^2(0,t;\boldsymbol{H})}
    \notag\\
&\quad \leq C\norma{\nabla \Delta \xi}_{L^2(0,t;\HHH)} +C \norma{\xi}_{L^2(0,t;{\Hx2})}
+C\norma{\boldsymbol{P}\bh}_{L^2(0,T;\boldsymbol{V}_\sigma^*)}^2,
\quad \forall\, t \in [0,T].
\end{align*}
Combining the above estimates on $G$ and its gradient, we arrive at the conclusion \eqref{eq:fre:G:H1}.

This completes the proof of Theorem \ref{THM:FRE}.
\end{proof}

\begin{remark}\rm
In fact, the proof of Theorem \ref{THM:FRE} implies the \Fre\ differentiability of $\Sol$ as a map $\bu \mapsto (\bv,\ph,\mu,\omega)$ between $\UR$ and a suitable Banach space. This property enables us to include in the cost functional \eqref{def:cost} additional terms of the form
\begin{align*}
	& \frac {\alpha_4}2 \intQ |\bv-\bv_Q|^2\,\mathrm{d}x\mathrm{d}t
	+\frac {\alpha_5}2 \iO |\bv(T)-\bv_\Omega|^2\,\mathrm{d}x
	+\frac {\alpha_6}2 \intQ |\mu-\mu_Q|^2\,\mathrm{d}x\mathrm{d}t
	+\frac {\alpha_7}2 \iO |\mu(T)-\mu_\Omega|^2\,\mathrm{d}x
	\\ & \quad
	+\frac {\alpha_8}2 \intQ |\omega-\omega_Q|^2\,\mathrm{d}x\mathrm{d}t
	+\frac {\alpha_9}2 \iO |\omega(T)-\omega_\Omega|^2\,\mathrm{d}x,
\end{align*}
for suitable target functions $\bv_Q$, $\mu_Q$, $\omega_Q$, $\bv_\Omega$, $\mu_\Omega$, $\omega_\Omega$, defined in $Q$ and $\Omega$, respectively,  and nonnegative coefficients $\alpha_i$, $i=4,...,9$. The entire theory for the corresponding optimal control problem could be developed {\it mutatis mutandis} without obstructions. However, since such terms are not particularly pertinent from a modeling perspective here, we opt to avoid including those terms for the sake of simplicity.
\end{remark}

\section{The Optimal Control Problem}
\setcounter{equation}{0}
\label{SEC:CONTROL}
Let us now study the optimal control problem \textbf{(CP)}. Our focus lies on addressing two fundamental questions: first, establishing the existence of at least one (globally) optimal control, and second, delineating first-order necessary optimality conditions.

\subsection{Existence of optimal controls}
\label{SUBSEC:EXMIN}

The existence of an optimal control can be achieved through an application of the direct method of calculus of variations. Given the standard nature of the technique, we proceed quickly, deferring to, e.g., \cite{CGSS6, PS} for further details, where a similar strategy has been employed for some comparable systems.

\begin{proof}[Proof of Theorem \ref{THM:EX:CONTROL}]
With the aid of the control-to-state mapping $\Sol$, we introduce the {\it reduced cost functional}
\begin{align}
	\label{def:Jred}
	\Jred(\bu):= {\cal J}(\bu,\Sol(\bu)),
	\quad \bu \in \mathcal{U},
\end{align}
so that the original control problem \textbf{(CP)} can be equivalently rephrased into the following minimization problem
$$
\min_{\bu \in \Uad} \Jred (\bu).
$$

First, we notice that the reduced cost functional $\Jred$ defined in \eqref{def:Jred} is nonnegative, whence it is bounded from below.
As a consequence, the infimum $\mathcal{J}^*=\inf_{\bu \in \Uad} \Jred (\bu)$ exists and it is finite. Moreover, there exists a minimizing sequence of controls $\{\bu_n\} \subset \Uad$ such that $\Jred (\bu_n)\to \mathcal{J}^*$ as $n\to \infty$.
 We denote the corresponding sequence of states as $\{\soln\}$ with $n\in\enne$. Then, up to a non-relabelled subsequence, there exists a limit $\bu^*$ such that, as $n\to\infty$,
\begin{align*}
	\bu_n \to \bu^* \quad \text{weakly star in} \,\,  \LLL^\infty(Q).
\end{align*}
Since $\Uad$ is convex and closed in $\mathcal{U}$, and thus weakly sequentially closed, we have $\bu^* \in \Uad$.
Exploiting the properties of the control-to-state operator $\Sol$ established in Theorem \ref{THM:RES}, we can prove the existence of a function $\ph^*$ such that, up to a non-relabeled subsequence, as $n\to\infty$, $\ph_n \to \ph^*$ weakly in the regularity space illustrated in Theorem \ref{THM:RES}.
Due to the sufficiently high regularity of $\ph_n$, it easily follows from compactness arguments that, up to a non-relabeled subsequence, as $n \to \infty$,
\begin{align*}
	\ph_n \to \ph^*
	\quad \text{strongly in} \,\,
	C^0(\ov Q).
\end{align*}
This enables us to infer that $\ph^*=\Sol(\bu^*)$. Therefore, since the sequence is minimizing, it readily follows from the weak lower semicontinuity of ${\cal J}$ and the above convergence results that
$$
{\cal J}(\bu^*,\Sol(\bu^*))
=\Jred(\bu^*)
\leq \liminf_{n\to \infty}\Jred(\bu_n)
=\mathcal{J}^*,
$$
namely, the pair $(\ph^*,\bu^*)$ yields an admissible optimal solution to the control problem \textbf{(CP)}.
\end{proof}

\subsection{First-order necessary optimality conditions}
\label{SUBSEC:FOC}

We now proceed to derive first-order necessary optimality conditions for locally optimal solutions (i.e., local minimizers of the cost functional).
First, we prove Theorem \ref{THM:VAR:INEQ:PREL} that  presents a necessary optimality condition in the form of a variational inequality.

\begin{proof}[Proof of Theorem \ref{THM:VAR:INEQ:PREL}]
Since the control-to-state operator $\Sol$ is Fr\'{e}chet
differentiable, so is the reduced cost functional $\Jred$ due to the chain rule. Invoking the convexity of $\Uad$, we find
\begin{align}
	\label{opt:abs}
(\Jred'(\widetilde{\bu}), \bu-\widetilde{\bu})\geq 0,\quad \forall\, \bu \in \Uad.
\end{align}
On the other hand, we infer from the chain rule that
$$
\Jred'(\bu)= \mathcal{J}'_{\mathcal{S}(\bu)}(\mathcal{S}(\bu), \bu)\circ D\mathcal{S}(\bu)+\mathcal{J}'_{\bu}(\mathcal{S}(\bu), \bu),
$$ where for every fixed $\bu\in \mathcal{U}$,
$\mathcal{J}'_{y}(y, \bu)$ is the Fr\'echet derivative of $\mathcal{J}(y,\bu)$ with respect to $y$ at $y\in \mathcal{Z}$, with $\cal Z$ being defined in \eqref{def:Z}, and for every fixed $y\in \mathcal{Z}$,
$\mathcal{J}'_{\bu}(y, \bu)$ denotes the Fr\'echet derivative with respect to $\bu$ at $\bu\in \mathcal{U}$. Using a straightforward computation and Theorem \ref{THM:FRE}, we obtain the variational inequality \eqref{var:ineq:prel}. This completes the proof.
\end{proof}

In what follows, we intend to express the optimality condition obtained in Theorem \ref{THM:VAR:INEQ:PREL} by means of the adjoint variables. For this purpose, we first establish the well-posedness of the adjoint system \Adj. The nonconstant mobility introduces several technical complexities that, while significant for the modeling, are not central to our current focus on the optimal control problem. To shorten the presentation and remain aligned with the objectives of our study, we opt to consider a constant mobility function for the rest of the analysis, taken equal to unity for simplicity. It is worth mentioning that analogous results can be obtained for the general case and the derivation of the adjoint system will be done considering a nonconstant mobility (see Remark~\ref{REM:MOB} below).

\begin{proposition}	\label{THM:ADJ}
Suppose that the assumptions in Theorem \ref{THM:EX:CONTROL} are satisfied and $m(\ph)\equiv 1$. Then, the adjoint system \Adj\ admits a unique weak solution $(\bva,\pha,\mua, \oma)$ on $[0,T]$ such that it enjoys the regularity properties
\begin{align*}
		\bva & \in  \L\infty {\Vs} \cap \L2 {\Ws}\cap \H1 {\Hs},
		\\
		\pha &\in  \L\infty H \cap \L2 {{ H^3(\Omega)}\cap \Hdn}\cap \H1 {(H^3(\Omega)\cap\Hdn)^*},
		\\
		\mua & \in \L2 V,\quad
		\oma  \in \L2 {V^*},
\end{align*}
	and fulfills the variational identities
\begin{align}
	\non
	& -\iO \dt \bva \cdot \bz \,\mathrm{d}x
	 - 2 \iO  \div(\nu(\ph)D\bva) \cdot \bz \,\mathrm{d}x
	- \iO (\bv\cdot \nabla)\bva \cdot \bz \,\mathrm{d}x
	+ \iO (\bva \cdot \nabla^\top)\bv\cdot \bz \,\mathrm{d}x
	\\ & \quad
	+ \iO\pha \nabla \ph \cdot \bz \,\mathrm{d}x
	=
	0 \quad \text{a.e. in $(0, T )$ and for every $\bz \in \Hs$},
    \label{we:adj:1}
    	\\
	& \non
	-\<\dt \pha, \xi>_{{H^3(\Omega)\cap \Hdn}}
	- \< \oma, \Delta \xi >_V
	+ 2 \iO \nu'(\ph)D\bv : \nabla \bva \,  \xi\,\mathrm{d}x
	- \iO (\nabla \pha \cdot \bv)\xi\,\mathrm{d}x
	\\
    & \qquad \non
	+ \iO (\nabla \mu \cdot \bva)\xi\,\mathrm{d}x
	+ \iO f''(\ph)\omega \mua \xi\,\mathrm{d}x
	+ \iO f'(\ph)\oma \xi\,\mathrm{d}x
	\\ & \quad 	
	=\alphaphQ \iO (\ph-\ph_Q)\xi\,\mathrm{d}x
	\quad \text{a.e. in $(0, T )$ and for every  $\xi \in H^3(\Omega)\cap \Hdn$},
    \label{we:adj:2}
    \\
    & \<\oma,\rho>_V
	= \int_\Omega \nabla \mua\cdot \nabla \rho\,\mathrm{d}x + \int_\Omega (f'(\ph)+\eta)\mua \rho\,\mathrm{d}x
	\quad \text{a.e. in $(0, T )$ and for every  $\rho \in V$}.
    \label{we:adj:3}
\end{align}
Besides, equation \eqref{eq:ad:4} is satisfied almost everywhere in $Q$, and the terminal conditions in \eqref{eq:ad:7} are fulfilled almost everywhere in $\Omega$, respectively.
\end{proposition}
\begin{remark}\rm
Similar consideration for the Lagrange multiplier related to the pressure $\pa$ can be made in the spirit of Remark \ref{RMK:pressure}.
Thanks to the Aubin--Lions--Simon theorem (see \cite{Simon}) and the regularity properties obtained in Proposition \ref{THM:ADJ}, we find $ \bva \in \C0{\Vs}$ and $\pha \in  \C0{H}$ so that the terminal data can be attained.
\end{remark}
\begin{proof}[Proof of Proposition \ref{THM:ADJ}]
Given our current focus on the constant mobility $m(\ph)\equiv1$, we obtain certain simplifications in the adjoint system \Adj, which motivate the form of the weak formulation proposed above.

The existence of a weak solution follows from a suitable Faedo--Galerkin method. We sketch it here for the convenience of the readers. Consider the family of eigenvalues $\{\lambda_i\}_{i=1}^{\infty}$ and the corresponding eigenfunctions $\{\bm{y}_{i}(x)\}_{i=1}^{\infty}$ with $\|\bm{y}_{i}\|=1$ of the Stokes operator $\bm{A}$.
It is well-known that $0<\lambda_1\leq \lambda_2 \leq \cdots \to +\infty$, and that $\{\bm{y}_{i}\}_{i=1}^{\infty}$ forms a complete orthonormal basis of $\Hs$, and is orthogonal in $\Vs$.
Next, we consider the family of eigenvalues $\{\ell_i\}_{i=1}^{\infty}$ and the corresponding eigenfunctions $\{z_{i}(x)\}_{i=1}^{\infty}$ with $\|z_{i}\|=1$ of the elliptic operator $I-\Delta$ subject to homogeneous Neumann boundary conditions.
Then $0<\ell_1<\ell_2 \leq \cdots \to +\infty$, and that $\{z_{i}\}_{i=1}^{\infty}$ forms a complete orthonormal basis of $H$, and is orthogonal in $V$.
For every integer $k\geq 1$, we denote the finite-dimensional subspace of $\Hs$ by
$\bm{Y}_{k}:=\textrm{span} \big\{\bm{y}_{1}(x) ,\cdots,\bm{y}_{k}(x)\big\}$. The orthogonal projection on $\bm{Y}_{k}$ with respect to the inner product in $\Hs$ is denoted by $\bm{P}_{\bm{Y}_{k}}$. Similarly, we denote the finite-dimensional subspace of $H$ by $Z_{k}:=\textrm{span} \big\{z_{1}(x) ,\cdots,z_{k}(x)\big\}$.
The orthogonal projection on $Z_{k}$ with respect to the inner product in $H$ is denoted by $\bm{P}_{Z_{k}}$.
We note that $\bigcup_{k=1}^\infty \bm{Y}_{k}$ is dense in $\Hs$, $\Vs$ and $\Ws$, while $\bigcup_{k=1}^\infty Z_{k}$ is dense in $H$, $V$, $\Hdn$ and $H^3(\Omega)\cap \Hdn$, respectively.

Let now $k\in \mathbb{Z}^+$. We look for finite dimensional approximate solutions in the form
$$
\bv^a_k=\sum_{i=1}^k a_{k,i}(t) \bm{y}_i(x),\quad
\varphi^a_k=\sum_{i=1}^k b_{k,i}(t)z_i(x),\quad
\mu^a_k=\sum_{i=1}^k c_{k,i}(t) z_i(x),\quad
\omega^a_k=\sum_{i=1}^k d_{k,i}(t) z_i(x),
$$
that satisfy
\begin{align}
	\non
	& -\iO \dt \bv^a_k \cdot \bz\,\mathrm{d}x
	 - 2 \iO  \div(\nu(\ph)D\bv^a_k) \cdot \bz \,\mathrm{d}x
	- \iO (\bv\cdot \nabla)\bv^a_k \cdot \bz \,\mathrm{d}x
	+ \iO (\bv^a_k \cdot \nabla^\top)\bv\cdot \bz \,\mathrm{d}x
	\\ & \quad
	+ \iO\varphi^a_k \nabla \ph \cdot \bz \,\mathrm{d}x
	=
	0 \quad \text{a.e. in $(0, T )$ and for every $\bz \in \bm{Y}_{k}$},
    \label{we:adj:1-ap}
    	\\
	& \non
	-\iO \dt \varphi^a_k \xi\,\mathrm{d}x
	- \iO \omega^a_k \Delta \xi \,\mathrm{d}x
	+ 2 \iO \nu'(\ph)D\bv : \nabla \bv^a_k \,  \xi\,\mathrm{d}x
	- \iO (\nabla \varphi^a_k \cdot \bv)\xi\,\mathrm{d}x
	\\
    & \qquad \non
	+ \iO (\nabla \mu \cdot \bv^a_k)\xi\,\mathrm{d}x
	+ \iO f''(\ph)\omega \mu^a_k \xi\,\mathrm{d}x
	+ \iO f'(\ph)\omega^a_k \xi\,\mathrm{d}x
	\\ & \quad 	
	=\alphaphQ \iO (\ph-\ph_Q)\xi\,\mathrm{d}x
	\quad \text{a.e. in $(0, T )$ and for every  $\xi \in Z_k$},
    \label{we:adj:2-ap}
    \\
    & \iO \omega^a_k \rho \,\mathrm{d}x
	= \int_\Omega \nabla \mu^a_k\cdot \nabla \rho\,\mathrm{d}x + \int_\Omega (f'(\ph)+\eta)\mu^a_k \rho\,\mathrm{d}x
	\quad \text{a.e. in $(0, T )$ and for every  $\rho \in Z_k$},
    \label{we:adj:3-ap}
\end{align}
and
\begin{align}
&\mu^a_k = -\Delta \varphi^a_k -\bm{P}_{Z_{k}}(\nabla  \ph \cdot \bv^a_k),
	\label{eq:ad:5-ap}\\
	& \bv^a_k(T)={\boldsymbol 0},
	\quad
	\varphi^a_k(T)=\alphaphO \bm{P}_{Z_{k}}(\ph(T)-\ph_\Omega).\label{eq:ad:7-ap}
\end{align}
%
Like before, it will be convenient to make the variable $\mu^a_k$ explicit in the definition of $\omega^a_k$, that is,
\begin{align}
	 \iO \omega^a_k \rho \,\mathrm{d}x
	& =
	 -\iO \nabla \Delta \varphi^a_k \cdot\nabla \rho\,\mathrm{d}x
	 - \iO \nabla (\bm{P}_{Z_{k}}(\nabla  \ph \cdot \bv^a_k))\cdot \nabla \rho\,\mathrm{d}x
     \notag\\
     &\quad - \int_\Omega (f'(\ph)+\eta)(\Delta \varphi^a_k +\bm{P}_{Z_{k}}(\nabla  \ph \cdot \bv^a_k)) \rho\,\mathrm{d}x
     \quad \text{a.e. in $(0, T )$ and for every  $\rho \in Z_k$}.
     \label{eq:ad:oma-ap}
\end{align}
From the classical Cauchy--Lipschitz
theorem for systems of ordinary differential equations, we can easily obtain the existence of a unique finite-dimensional approximate solution $(\bv^a_k,\varphi^a_k,\mu^a_k,\omega^a_k)$ in the time interval $[0,T]$.

Next, we derive estimates that are uniform with respect to the parameter $k$. For this purpose, we test \eqref{we:adj:1-ap} with $\bv^a_k + \bm{A} \bv^a_k$, \eqref{we:adj:2-ap} with $\varphi^a_k$, and \eqref{eq:ad:oma-ap} with $\Delta \varphi^a_k$. Upon adding the resulting identities and rearranging terms, we infer that
\begin{align}
	& - \frac 12 \frac{\mathrm{d}}{\mathrm{d}t} \big(\norma{\bv^a_k}_{\Vs}^2   +\norma{\varphi^a_k}^2\big)
        + 2\iO \nu(\varphi)|D \bv^a_k|^2 \,\mathrm{d}x
	+ \iO   |\nabla \Delta \varphi^a_k|^2\,\mathrm{d}x
	\notag \\
    & \quad = 2 \iO  \div(\nu(\ph)D\bv^a_k) \cdot \bm{A} \bv^a_k \,\mathrm{d}x
     + \iO (\bv\cdot \nabla)\bv^a_k \cdot (\bv^a_k + \bm{A} \bv^a_k) \,\mathrm{d}x
     \notag \\
     &\qquad
	- \iO (\bv^a_k \cdot \nabla^\top)\bv\cdot (\bv^a_k + \bm{A} \bv^a_k) \,\mathrm{d}x
	- \iO \varphi^a_k \nabla \ph \cdot (\bv^a_k + \bm{A} \bv^a_k) \,\mathrm{d}x
	\notag \\
       & \qquad
	- 2 \iO \nu'(\ph)(D\bv : \nabla \bv^a_k)\varphi^a_k \,\mathrm{d}x
	+  \iO (\nabla\varphi^a_k \cdot \bv) \varphi^a_k \,\mathrm{d}x
    \notag \\
    &\qquad
	- \iO (\nabla \mu \cdot \bv^a_k) \varphi^a_k \,\mathrm{d}x
	- \iO f''(\ph)\omega \mu^a_k \varphi^a_k \,\mathrm{d}x
    - \iO f'(\ph)\omega^a_k \varphi^a_k\,\mathrm{d}x
    \notag \\
    &\qquad 	
	+ \alphaphQ \iO (\ph-\ph_Q)\varphi^a_k\,\mathrm{d}x
	- \iO \nabla (\bm{P}_{Z_{k}}(\nabla  \ph \cdot \bv^a_k)) \cdot \nabla \Delta \varphi^a_k \,\mathrm{d}x
    \notag \\
    &\qquad
	- \iO   ( f'(\ph)	+\eta)( \Delta \varphi^a_k  + \bm{P}_{Z_{k}}(\nabla  \ph \cdot \bv^a_k)) \Delta \varphi^a_k \,\mathrm{d}x	
	\notag \\
    & \quad
	 =:\sum_{j=1}^{12}I_j.
     \label{ap-diff}
\end{align}
We then proceed to estimate the terms on the \rhs\ of \eqref{ap-diff}. Using the identity
\begin{align*}
	2\div ( \nu(\ph) D\bv^a_k)
= 2\nu'(\ph)  D\bv^a_k  \nabla \ph +   \nu(\varphi)\Delta \bv^a_k,
\end{align*}
and \ref{ass:1:visco}, we can apply the same argument as in \cite[Section 4]{WY} to conclude
\begin{align*}
I_1\leq -\frac{\nu_*}{2}\|\bm{A}\bv^a_k\|^2+ C \|\nabla \bv^a_k\|^2.
\end{align*}
%
%
Next, using the \eqref{control:reg} for $(\bv,\varphi, \mu,\omega)$, the Sobolev embedding theorem, Lemma \ref{sto}, the H\"older, Young, Poincar\'e, and Korn inequalities, for a constant $\delta\in (0,1)$ yet to be chosen, we have
\begin{align*}
 I_2 & \leq
	 \norma{\bv }_4 \norma{\nabla\bv^a_k} \norma{\bv^a_k}_4
	 +\norma{\bv }_4 \norma{\nabla\bv^a_k}_4  \norma{ \bm{A}\bv^a_k}
	\\ &
	\leq
	 C \norma{\bv }_{\Vs} \norma{\bv^a_k}_{\Vs}^2
	 + C \norma{\bv }_{\Vs} \norma{\nabla\bv^a_k}^{\frac{1}{2}} \norma{ \bm{A} \bv^a_k}^{\frac{3}{2}}
	\\ &
	\leq C \norma{\bv^a_k}^2_{\Vs}
        + \delta \norma{ \bm{A} \bv^a_k}^2	
	+ \cd \norma{\bv }_{\Vs}^4 \norma{\nabla\bv^a_k}^2
    \\ & \leq
	\delta \norma{ \bm{A} \bv^a_k}^2	
	+ \cd \norma{\bv^a_k}^2_{\Vs},
\end{align*}
\begin{align*}
I_3 & \leq
	 \norma{\bv^a_k }_4 \norma{\nabla\bv}_4 (\norma{\bv^a_k} +  \norma{ \bm{A} \bv^a_k})
	\leq
	\delta \norma{\bm{A} \bv^a_k}^2
	+ \cd (1+\norma{\bv}_{\HHH^2(\Omega)}^2)\norma{ \bv^a_k}^2_{\Vs},
\end{align*}
\begin{align*}
I_4 & \leq
	  	\norma{\varphi^a_k } \norma{\nabla \ph }_\infty (\norma{\bv^a_k} + \norma{ \bm{A} \bv^a_k})
	  	\leq
	\delta \norma{ \bm{A} \bv^a_k}^2
	+ \cd \norma{\varphi^a_k}^2
	+ C \norma{\bv^a_k}^2,
\end{align*}
\begin{align*}
I_5 & \leq
	  2\norma{\nu'(\ph)}_\infty \norma{D\bv}_4 \norma{\nabla \bv^a_k} \norma{\varphi^a_k}_4
      \\
      & \leq C \norma{\bv}_{\HHH^2(\Omega)} \norma{\bv^a_k}_{\Vs}\norma{\varphi^a_k}_{V}
      \\
      & \leq 	
	  \delta \norma{\nabla \Delta \varphi^a_k}^2
        + \cd \norma{\varphi^a_k}^2
	  + C \norma{\bv}^2_{\HHH^2(\Omega)}   \norma{\bv^a_k}^2_{\Vs},
\end{align*}
\begin{align*}
I_6 & \leq
	  \norma{\nabla \varphi^a_k }_4 \norma{\bv }_4 \norma{\varphi^a_k}
       \leq
	  C \norma{\bv }_{\Vs} \norma{\varphi^a_k}_{\Hx2}^2
	 \leq \delta \norma{\nabla \Delta \varphi^a_k}^2 + \cd \norma{\varphi^a_k}^2,
\end{align*}
\begin{align*}
I_7 & \leq 	
	\norma{\nabla \mu }_4\norma{\bv^a_k}_4  \norma{\varphi^a_k}
      \leq  C \norma{\varphi^a_k}^2
	  + C  \norma{ \mu }_{\Hx2}^2\norma{\bv^a_k}_{\Vs}^2,
\end{align*}
\begin{align*}
 I_8 & \leq 	
	  \norma{f''(\ph)}_\infty \norma{\omega}_\infty(\norma{\Delta \varphi^a_k} + \norma{\bm{P}_{Z_k}(\nabla \ph\cdot \bv^a_k)}) \norma{\varphi^a_k}
      \\
       & \leq
	  C \norma{\Delta \varphi^a_k} \norma{\varphi^a_k} + C\norma{\nabla \varphi}_\infty\norma{\bv^a_k} \norma{\varphi^a_k}
      \\
      & \leq \delta \norma{\nabla \Delta \varphi^a_k}^2
 	 + \cd (\norma{\varphi^a_k}^2
 	 + \norma{\bv^a_k}^2).
\end{align*}
To estimate the term $I_{9}$, we need to work a little bit more using \eqref{eq:ad:oma-ap}.
To this end, we observe that
\begin{align*}
	I_9 &
	=	\iO \nabla  \Delta \varphi^a_k
    \cdot\nabla ( \bm{P}_{Z_k}(f'(\ph) \varphi^a_k))\,\mathrm{d}x
	 + \iO \nabla (\bm{P}_{Z_k}(\nabla \ph  \cdot \bv^a_k))\cdot \nabla ( \bm{P}_{Z_k}(f'(\ph) \varphi^a_k))\,\mathrm{d}x
	 \\ & \quad
	 + \iO
	 ( f'(\ph)	+\eta)(\Delta \varphi^a_k + \bm{P}_{Z_k}(\nabla  \ph \cdot \bv^a_k)) ( \bm{P}_{Z_k}(f'(\ph) \varphi^a_k))\,\mathrm{d}x
     \\
     & =: I_{9}^{(1)}+I_{9}^{(2)}+I_{9}^{(3)},
\end{align*}
with the following estimates
\begin{align*}
	I_{9}^{(1)}
	& \leq
	\norma{\nabla \Delta \varphi^a_k}\norma{ f'(\varphi)\varphi^a_k}_V
	  \leq \delta \norma{\nabla \Delta \varphi^a_k}^2
	  + \cd \norma{\varphi^a_k}^2,
\end{align*}
\begin{align*}
	I_{9}^{(2)}
	& \leq
	  \|\nabla \varphi\cdot \bv^a_k\|_V\norma{ f'(\varphi)\varphi^a_k}_V
       \leq
	  \delta \norma{\nabla \Delta \varphi^a_k}^2
	 + \cd \norma{\varphi^a_k}^2
 	 + C \norma{\bv^a_k}^2_{\Vs},
\end{align*}
\begin{align*}
	I_{9}^{(3)}
	& \leq 	
	\norma{f'(\ph)	+\eta}_\infty(\norma{\Delta \varphi^a_k  } + \norma{ \bm{P}_{Z_k}(\nabla  \ph \cdot \bv^a_k)})  \norma{  f'(\ph)}_\infty\norma{ \varphi^a_k}
	\\
        & \leq C \norma{\Delta \varphi^a_k}^2
	+ C  \norma{\varphi^a_k}^2 + \norma{\bv^a_k}^2
        \\
        & \leq \delta \norma{\nabla \Delta \varphi^a_k}^2
	+ \cd \norma{ \varphi^a_k}^2
	+ C  \norma{\bv^a_k}^2.
\end{align*}
The term $I_{10}$ can be easily handled, that is,
\begin{align*}
 I_{10} & \leq 	
	  \alphaphQ\norma{\ph-\ph_Q}(\norma{\varphi^a_k} + \norma{\Delta \varphi^a_k})
      \\
	  &\leq  \delta \norma{\nabla \Delta\varphi^a_k}^2
 	 + \cd  \norma{ \varphi^a_k}^2
 	 + \cd  (\norma{ \ph}^2+ \norma{ \ph_Q}^2).
\end{align*}
Next, we have
\begin{align*}
I_{11} & \leq \norma{\nabla (\nabla \varphi\cdot\bv^a_k)}\norma{\nabla \Delta\varphi^a_k}
\\
&\leq (\norma{\nabla^2\varphi}_4\|\bv^a_k\|_4 +   \norma{\nabla\varphi}_\infty\|\nabla \bv^a_k\|) \norma{\nabla \Delta\varphi^a_k}
\\
&\leq C\|\bv^a_k\|_{\Vs}\norma{\nabla \Delta\varphi^a_k}
\\
&\leq \delta  \norma{\nabla \Delta\varphi^a_k}^2 + \cd \|\bv^a_k\|_{\Vs}^2,
\end{align*}
and finally,
\begin{align*}
    I_{12} & \leq 	
	  \norma{f'(\ph)+\eta}_{\infty}( \norma{ \Delta \varphi^a_k } + \norma{ \bm{P}_{Z_k}(\nabla  \ph \cdot \bv^a_k)}) \norma{\Delta\varphi^a_k}
      \\
      & \leq C \norma{ \Delta\varphi^a_k}^2
       + C \norma{\bv^a_k}^2
	\\ &
	  \leq  \delta \norma{\nabla \Delta \varphi^a_k}^2
	  + \cd \norma{\varphi^a_k}^2
	   + \cd \norma{\bv^a_k}^2.
\end{align*}
At this stage, we select $\delta\in (0,1)$ sufficiently small and then integrate \eqref{ap-diff} over time, using \eqref{control:reg}, \eqref{eq:ad:7} and assumptions \ref{ass:1:visco}, \ref{ass:6:cost}, we can apply (backward) Gronwall's lemma and the elliptic regularity theory to conclude that
\begin{align*}
	& \norma{\bv^a_k}_{\L\infty {\Vs} \cap \L2 {\Ws}} +\norma{\varphi^a_k}_{\L\infty {H} \cap \L2 {\Hx3}} \leq C,
\end{align*}
where the constant $C>0$ is independent of the approximating parameter $k$. With this, it readily follows, arguing by comparison in equations \eqref{eq:ad:5-ap} that
\begin{align*}
	\norma{\mu^a_k}_{\L2 V}\leq C,
\end{align*}
whence, due to the elliptic regularity theory and \eqref{we:adj:3-ap}, it holds
\begin{align*}
	\norma{\omega^a_k}_{\L2 {V^*}}\leq C.
\end{align*}
Finally, using \eqref{we:adj:1-ap} and \eqref{we:adj:2-ap}, we infer that
\begin{align*}
		& \norma{ \dt \bv^a_k}_{\L2 \Hs}\leq C,
		\quad
		\norma{ \dt \varphi^a_k}_{\L2 {(H^3(\Omega)\cap\Hdn)^*}}\leq C.
\end{align*}
From the above uniform estimates for approximate solutions, it is standard to apply the compactness method to extract a convergent subsequence $\{(\bv^a_k,\varphi^a_k,\mu^a_k,\omega^a_k)\}$ (not relabeled for simplicity) such that after passing $k\to +\infty$, the corresponding limit $\{(\bv^a,\varphi^a,\mu^a,\omega^a)\}$ yields a weak solution to the adjoint system \Adj\ on $[0,T]$. The details are omitted here.

Given the linearity of the problem \Adj\ and the regularity of the solution $(\bva,\pha,\mua, \oma)$, the uniqueness readily follows.
\end{proof}

\begin{remark}\label{REM:MOB}\rm
The methodology of the proof for the general mobility case remains essentially the same; however, numerous additional lower-order terms appear that have to be treated. Specifically, we can still express $\oma$ utilizing the equation for $\mua$, yielding:
\begin{align*}
	 \oma &=
	m(\ph) \Delta^2 \pha
	+ \Delta (\nabla \ph \cdot \bva)
	 -(f'(\ph)+\eta)
	(m(\ph) \Delta \pha
	+ \nabla \ph \cdot \bva)
      + L(\pha;\ph),
\end{align*}
where
\begin{align*}
	L(\pha;\ph)  &:= \Delta(m'(\varphi)\nabla \varphi\cdot\pha) + m''(\varphi)|\nabla \varphi|^2\Delta \pha + m'(\varphi)\Delta\varphi\Delta\pha+ 2m'(\varphi)\nabla \varphi\cdot\nabla \Delta \pha\\
&\quad \ -(f'(\ph)+\eta)
	m'(\ph) \nabla\ph \cdot \nabla \pha.
\end{align*}
If the mobility $m$ is constant, then $L(\pha;\ph) =0$, which reduces to the case we have treated. The details for the general case are left to the interested readers.
\end{remark}

We are now in a position to prove Theorem \ref{THM:VAR:INEQ:FINAL}.
\begin{proof}[Proof of Theorem \ref{THM:VAR:INEQ:FINAL}]
The conclusion follows from a combination of the variational inequality \eqref{var:ineq:prel} along with the adjoint system \Adj.
The following computations also provide a (formal) way to derive the adjoint system corresponding to the state system \Sys.
For completeness, we proceed with  the derivation for the general case when the mobility $m$ is phase dependent instead of constant, even if the corresponding analysis was carried out for the simplified case $m(\varphi)\equiv 1$ (cf. Proposition \ref{THM:ADJ} and Remark \ref{REM:MOB}).

Let $\opt\in\Uad$ be an optimal control with the associate state $\solopt$ and the adjoint state $\soladj$, which are the corresponding solutions to the systems \Sys\ and \Adj, respectively. For arbitrary $\bu\in\UR$, we denote $\bh = \bu - \opt$ and consider $\sollin$ the associated solution to the linearized system \eqref{eq:lin:1}--\eqref{eq:lin:7}. Then, we multiply \eqref{eq:lin:1} by $\bva$,
\eqref{eq:lin:2} by $p^a$,
\eqref{eq:lin:3} by $\pha$,
\eqref{eq:lin:4} by $\mua$,
and
\eqref{eq:lin:5} by $\oma$. Adding the resulting equalities, integrating over $(0,T)$ and $\Omega$, using integration by parts and also the initial conditions \eqref{eq:lin:7} as well as the terminal conditions \eqref{eq:ad:7}, we infer that
\begin{align}
	0 & =
	\ioT\<	\dt \bw , \bva>_{\Vs} \,\mathrm{d}t
    \notag \\
	& \quad +\intQ ((\bw \cdot \nabla ) \widetilde{\bv}
	+ (\widetilde{\bv} \cdot \nabla ) \bw
	- \div (2 \nu' (\widetilde{\ph}) \psi D \widetilde{\bv})
	- \div (2 \nu (\widetilde{\ph}) D \bw)) \cdot \bva \,\mathrm{d}x\mathrm{d}t
	\notag \\
    & \quad
	+ \intQ( \nabla q	
	- \theta \nabla \widetilde{\ph}
	- \widetilde{\mu} \nabla \psi -\bh)
    \cdot \bva\,\mathrm{d}x\mathrm{d}t
    + \intQ(\div \bw)\, p^a \,\mathrm{d}x\mathrm{d}t
	\notag \\
    & \quad
	+ \ioT \<\dt \psi, \pha>_{\Hdn}\,\mathrm{d}t
    \notag \\
    &\quad
	+ \intQ ( \bw \cdot \nabla \widetilde{\ph}
	+ \widetilde{\bv} \cdot \nabla \psi
	- \div (m'(\widetilde{\ph})\psi\nabla \widetilde{\mu})
	- \div (m(\widetilde{\ph})\nabla \theta)) \pha\,\mathrm{d}x\mathrm{d}t
	\notag \\
    & \quad
	+ \intQ (-\theta - \Delta \w
	+ f''(\widetilde{\ph})\psi \widetilde{\omega}
	+ f'(\widetilde{\ph}) \w
	+ \eta \w) \mua\,\mathrm{d}x\mathrm{d}t
	\notag \\
    & \quad
	+ \intQ (-\w - \Delta \psi
	+f'(\widetilde{\ph})\psi) \oma\,\mathrm{d}x\mathrm{d}t
	\notag \\
    & =
	- \intQ \bh \cdot \bva\,\mathrm{d}x\mathrm{d}t
	- \ioT\< \dt \bva , \bw>_{\Vs}\, \mathrm{d}t
	\notag \\
    & \quad
	+ \intQ (- \div (2 \nu(\widetilde{\ph})D\bva)
	- (\widetilde{\bv}\cdot \nabla)\bva
	+ (\bva \cdot \nabla^\top)\widetilde{\bv}
	+ \pha \nabla \widetilde{\ph}
    - \nabla p^a)\cdot \bw\,\mathrm{d}x\mathrm{d}t
	\notag \\
    & \quad
	-\ioT \<\dt \pha,\psi>_{\Hdn}\, \mathrm{d}t
	+ \iO  \psi(T)\pha(T) \,\mathrm{d}x
	\notag \\
    & \quad
	+ \intQ (- \Delta \oma
	+ 2 \nu'(\widetilde{\ph})D\widetilde{\bv} : \nabla \bva
	+ \nabla \widetilde{\mu} \cdot  \bva
	- \nabla \pha \cdot \widetilde{\bv}
	+ m'(\widetilde{\ph})\nabla \widetilde{\mu} \cdot \nabla \pha
	)\psi\,\mathrm{d}x\mathrm{d}t
	\notag \\
    & \quad
	+ \intQ (f''(\widetilde{\ph})\widetilde{\omega} \mua
	+ f'(\widetilde{\ph})\oma) \psi \,\mathrm{d}x\mathrm{d}t
	- \intQ (\mua + \div (m(\widetilde{\ph}) \nabla \pha)
     + \nabla \widetilde{\ph} \cdot \bva) \theta \,\mathrm{d}x\mathrm{d}t
	\notag \\
    & \quad
	+ \intQ (-\oma - \Delta \mua
    + f'(\widetilde{\ph})\mua
    + \eta \mua) \w \,\mathrm{d}x\mathrm{d}t
	\notag \\
	&
	= - \intQ \bh \cdot \bva \,\mathrm{d}x\mathrm{d}t
	+ \alphaphQ \intQ (\widetilde{\ph} - \ph_Q)\psi \,\mathrm{d}x\mathrm{d}t
	+ \alphaphO \iO (\widetilde{\ph}(T) - \ph_\Omega)\psi(T)\,\mathrm{d}x.
\label{adj-eq}
\end{align}
In the derivation of \eqref{adj-eq}, we have used the
integration-by-parts formula for
functions belonging to $\H1 \Vsp \cap \L2 \Vs$ and $\H1 {\Hdn^*} \cap \L2 \Hdn$,
and employ the definition of the adjoint system \Adj. We also point out that the following identities have been used
\begin{align*}
	(\bw \cdot \nabla ) \widetilde{\bv} \cdot \bva
	&
	= \sum_{i,j=1}^{2} \frac{\partial \widetilde{v}_i}{\partial {x_j}}w_j v_i^a
	= \sum_{i,j=1}^{2} \frac{\partial \widetilde{v}_j}{\partial {x_i}} v_j^a w_i
	= (\bva \cdot \nabla^\top ) \widetilde{\bv} \cdot \bw,
	\\
	(\widetilde{\bv} \cdot \nabla ) \bw \cdot \bva
	& = \sum_{i,j=1}^{2} \frac{\partial w_i}{\partial {x_j}}\widetilde{v}_j v_i^a
	= - \sum_{i,j=1}^{2} \frac{\partial v_i^a}{\partial {x_j}}\widetilde{v}_j w_i
	= - (\widetilde{\bv} \cdot \nabla ) \bva \cdot \bw,
\end{align*}
with natural notation for the vectors' components.

Combining the variational inequality \eqref{var:ineq:prel} and the identity \eqref{adj-eq}, we easily obtain the conclusion \eqref{var:ineq}.
The last part concerning the projection formula is a direct consequence of Hilbert's projection theorem.
The proof of Theorem \ref{THM:VAR:INEQ:FINAL} is complete.
\end{proof}

Let us conclude with a brief remark on numerical perspectives. By comparing the theoretical condition \eqref{opt:abs} with the explicit expression \eqref{var:ineq}, it follows, via Riesz’s representation theorem, that the gradient of the reduced cost functional $\Jred $ can be identified as
$$\nabla \Jred (\opt) = 	\alphacon \opt  + \bva.$$
This allows the optimal control problem to be viewed as a constrained minimization with a known gradient, naturally suggesting gradient-based methods, such as the projected conjugate gradient scheme. Nevertheless, the problem is challenging due to its strong nonlinearity, and its numerical implementation is therefore left as a particularly ambitious topic for future work. We refer to \cite{Dede,Hin} for further insights and computational perspectives in the context of optimal control problems.

\section*{Declarations}
\noindent
\textbf{Conflict of interest.} The authors have no competing interests to declare that are relevant to the content of this article.
\smallskip
\\
\noindent
\textbf{Data availability.} Data sharing not applicable to this article as no datasets were generated or analyzed during the current study.
\smallskip
\\
\noindent
\textbf{Ethical approval.} This research does not involve humans and/or animals.
\smallskip
\\
\noindent
\textbf{Fundings.} AS gratefully acknowledges partial
support by MUR, grant Dipartimento di Eccellenza 2023--2027, from the GNAMPA (Gruppo Nazionale per l'Analisi Matematica, la Probabilit\`a e le loro Applicazioni) of INdAM (Isti\-tuto Nazionale di Alta Matematica)
project CUP E5324001950001,
and from the Alexander von Humboldt Foundation.
\smallskip
\\
\noindent
\textbf{Acknowledgments.}
The authors thank the anonymous referee for helpful suggestions that allowed them to improve the article. AS acknowledges his affiliation to the GNAMPA (Gruppo Nazionale per l'Analisi Matematica, la Probabilit\`a e le loro Applicazioni) of INdAM (Isti\-tuto Nazionale di Alta Matematica). HW is a member of the Key Laboratory of Mathematics for Nonlinear Sciences (Fudan University), Ministry of Education of China.


\vspace{3truemm}

\Begin{thebibliography}{10}
\footnotesize
\itemsep=-3pt

\bibitem{Aland2017}
S. Aland,
Phase field models for two-phase flow with surfactants and biomembranes,
in \textit{Transport Processes at Fluidic Interfaces}, Adv. Math. Fluid Mech., Birkhäuser/Springer, Cham, (2017), 271--290.

\bibitem{AlandEgererLowengrubVoigt2014}
S. Aland, S. Egerer, J. Lowengrub, and A. Voigt,
Diffuse interface models of locally inextensible vesicles in a viscous fluid,
\textit{J. Comput. Phys.}, \textbf{277} (2014), 32--47.

\bibitem{Bo21}
F. Boyer and P. Fabrie,
\emph{Mathematical Tools for the Study of the Incompressible Navier--Stokes Equations and Related Models},
Applied Mathematical Sciences, Springer New York, 2013.

\bibitem{BCMS}
M. Burger, S. Y. Chu, P. A. Markowich, and C. B. Sch\"{o}nlieb,
The Willmore functional and instabilities in the Cahn--Hilliard equation,
{\it Commun. Math. Sci.}, {\bf 6} (2008), 309--329.

\bibitem{CampeloHernandezMachado2006}
F. Campelo and A. Hern\'{a}ndez-Machado,
Dynamic model and stationary shapes of fluid vesicle,
\textit{Eur. Phys. J. E}, \textbf{20} (2006), 37--45.

\bibitem{Canham1970}
P. Canham,
The minimum energy of bending as a possible explanation of the bioconcave shape of the human red blood cell,
\textit{J. Theoret. Biol.}, \textbf{26} (1970), 61--81.

\bibitem{ChengWangWiseYuan2020}
K. Cheng, C. Wang, S. Wise, and Z. Yuan,
Global-in-time Gevrey regularity solutions for the functionalized Cahn--Hilliard equation,
\textit{Discrete Contin. Dyn. Syst. Ser. S}, \textbf{13} (2020), 2211--2229.

\bibitem{ClimentEzquerraGuillenGonzalez2016}
B. Climent-Ezquerra and F. Guillén-González,
Long-time behavior of a Cahn--Hilliard--Navier--Stokes vesicle-fluid interaction model,
\textit{Trends in Differential Equations and Applications}, SEMA SIMAI Springer Ser., Springer, \textbf{8} (2016), 125--145.

\bibitem{ClimentEzquerraGuillenGonzalez2019}
B. Climent-Ezquerra and F. Guill\'{e}n-Gonz\'{a}lez,
Convergence to equilibrium of global weak solutions for a Cahn--Hilliard--Navier--Stokes vesicle model,
\textit{Z. Angew. Math. Phys.}, \textbf{70} (2019), Paper No. 125, 27 pp.

\bibitem{CGhyp}
P. Colli and G. Gilardi,
Hyperbolic relaxation of a sixth-order Cahn--Hilliard equation,
\textit{Discrete Contin. Dyn. Syst. Ser. S}, early access (2025), DOI: 10.3934/dcdss.2025106

\bibitem{ColliLaurencot2011}
P. Colli and P. Lauren\c{c}ot,
A phase-field approximation of the Willmore flow with volume constraints,
\textit{Interfaces Free Bound.}, \textbf{13} (2011), 341--351.

\bibitem{ColliLaurencot2012}
P. Colli and P. Lauren\c{c}ot,
A phase-field approximation of the Willmore flow with volume and area constraints,
\textit{SIAM J. Math. Anal.}, \textbf{44} (2012), 3734--3754.

\bibitem{CGSS6}
P. Colli, G. Gilardi, A. Signori, and J. Sprekels,
Curvature effects in pattern formation: well-posedness and optimal control of a sixth-order Cahn--Hilliard equation,
{\it SIAM J. Math. Anal.}, {\bf 56} (2024), 4928--4969.

\bibitem{ColliGilardiSprekels2018}
P. Colli, G. Gilardi, and J. Sprekels,
Optimal velocity control of a viscous Cahn--Hilliard system with convection and dynamic boundary conditions,
\textit{SIAM J. Control Optim.}, \textbf{56} (2018), 1665--1691.

\bibitem{CGSS8}
P. Colli, G. Gilardi, A. Signori, and J. Sprekels,
On Brinkman flows with curvature-induced phase separation in binary mixtures
Preprint arXiv:2509.20282 [math.AP], (2025), 1--29.

\bibitem{DaiLiuLuongPromislow2021}
S. Dai, Q. Liu, T. Luong, and K. Promislow,
On nonnegative solutions for the functionalized Cahn--Hilliard equation with degenerate mobility,
\textit{Results Appl. Math.}, \textbf{12} (2021), 100195, 13 pp.

\bibitem{DaiLiuPromislow2021}
S. Dai, Q. Liu, and K. Promislow,
Weak solutions for the functionalized Cahn--Hilliard equation with degenerate mobility,
\textit{Appl. Anal.}, \textbf{100} (2021), 1--16.

\bibitem{DaiPromislow2013}
S. Dai and K. Promislow,
Geometric evolution of bilayers under the functionalized Cahn--Hilliard equation,
\textit{Proc. Roy. Soc. A}, \textbf{469} (2013), 20120505, 20 pp.

\bibitem{DaiPromislow2015}
S. Dai and K. Promislow,
Competitive geometric evolution of amphiphilic interfaces,
\textit{SIAM J. Math. Anal.}, \textbf{47} (2015), 347--380.

\bibitem{Dede}
L. Ded\`{e},
Optimal flow control for Navier--Stokes equations: drag minimization,
\textit{Int. J. Numer. Method Fluid.}, \textbf{55} (2007), 347--366.

\bibitem{DHPW}
A. Doelman, G. Hayrapetyan, K. Promislow and B. Wetton,
Meander and pearling of single-curvature bilayer interfaces in the functionalized Cahn--Hilliard equation,
{\it SIAM J. Math. Anal.} {\bf 46} (2014), 3640--3677.

\bibitem{DuLiLiu2007}
Q. Du, M.-L. Li, and C. Liu,
Analysis of a phase field Navier--Stokes vesicle-fluid interaction model,
\textit{Discrete Contin. Dyn. Syst. Ser. B}, \textbf{8} (2007), 539--556.

\bibitem{DuLiuRyhamWang2005}
Q. Du, C. Liu, R. Ryham, and X.-Q. Wang,
A phase field formulation of the Willmore problem,
\textit{Nonlinearity}, \textbf{18} (2005), 1249--1267.

\bibitem{DuLiuRyhamWang2005b}
Q. Du, C. Liu, R. Ryham, and X.-Q. Wang,
Phase field modeling of the spontaneous curvature effect in cell membranes,
\textit{Commun. Pure Appl. Anal.}, \textbf{4} (2005), 537--548.

\bibitem{DuLiuRyhamWang2009}
Q. Du, C. Liu, R. Ryham, and X.-Q. Wang,
Energetic variational approaches in modeling vesicle and fluid interactions,
\textit{Phys. D}, \textbf{238} (2009), 923--930.

\bibitem{DuLiuWang2004}
Q. Du, C. Liu, and X.-Q. Wang,
A phase field approach in the numerical study of the elastic bending energy for vesicle membranes,
\textit{J. Comput. Phys.}, \textbf{198} (2004), 450--468.

\bibitem{DuLiuWang2006}
Q. Du, C. Liu, and X.-Q. Wang,
Simulating the deformation of vesicle membranes under elastic bending energy in three dimensions,
\textit{J. Comput. Phys.}, \textbf{212} (2006), 757--777.

\bibitem{EntringerBoldrini2015}
A. P. Entringer and J. L. Boldrini,
A phase field $\alpha$-Navier--Stokes vesicle-fluid interaction model: Existence and uniqueness of solutions,
\textit{Discrete Contin. Dyn. Syst. Ser. B}, \textbf{20} (2015), 397--422.

\bibitem{FRS}
S. Frigeri, E. Rocca, and J. Sprekels,
Optimal distributed control of a nonlocal Cahn--Hilliard/Navier--Stokes system in two dimensions.
{\it SIAM J. Control Optim.}, {\bf 54} (2016), 221--250.

\bibitem{GavishHayrapetyanPromislowYang2011}
N. Gavish, G. Hayrapetyan, K. Promislow, and L. Yang,
Curvature driven flow of bi-layer interfaces,
\textit{Phys. D}, \textbf{240} (2011), 675--693.

\bibitem{GavishJonesXuChristliebPromislow2012}
N. Gavish, J. Jones, Z. Xu, A. Christlieb, and K. Promislow,
Variational models of network formation and ion transport: Applications to perfluorosulfonate ionomer membranes,
\textit{Polymers}, \textbf{4} (2012), 630--655.

\bibitem{GigaKirshteinLiu2018}
M.-H. Giga, A. Kirshtein, and C. Liu,
Variational modeling and complex fluids,
in \textit{Handbook of Mathematical Analysis in Mechanics of Viscous Fluids}, Springer, Cham, (2018), 73--113.

\bibitem{GompperSchick1990}
G. Gompper and M. Schick,
Correlation between structural and interfacial properties of amphiphilic systems,
\textit{Phys. Rev. Lett.}, \textbf{65} (1990), 1116--1119.

\bibitem{GW1}
M. Grasselli and H. Wu,
Well-posedness and long-time behavior for the modified phase-field crystal equation,
{\it Math. Models Methods Appl. Sci.}, {\bf 24} (2014), 2743--2783.

\bibitem{GW2}
M. Grasselli and H. Wu,
Robust exponential attractors for the modified phase-field crystal equation,
{\it Discrete Contin. Dyn. Syst.}, {\bf 35} (2015), 2539--2564.

\bibitem{Helfrich1973}
W. Helfrich,
Elastic properties of lipid bilayers: Theory and possible experiments,
\textit{Z. Naturforsch. C}, \textbf{28} (1973), 693--703.

\bibitem{Hin}
M. Hinterm\"uller, M. Hinze, C. Kahle, and T. Keil,
A goal-oriented dual-weighted adaptive finite element approach for the optimal control of a nonsmooth Cahn--Hilliard--Navier--Stokes system,
{\it Optim. Eng.}, {\bf 19} (2018), 629--662.

\bibitem{KNR}
M. Korzec, P. Nayar, and P. Rybka,
Global weak solutions to a sixth order Cahn--Hilliard type equation,
{\it SIAM J. Math. Anal.}, {\bf 44} (2012), 3369--3387.

\bibitem{KR}
M. Korzec and P. Rybka,
On a higher order convective Cahn--Hilliard-type equation,
{\it SIAM J. Appl. Math.}, {\bf 72} (2012), 1343--1360.

\bibitem{LiuTakahashiTucsnak2012}
Y.-N. Liu, T. Takahashi, and M. Tucsnak,
Strong solution for a phase-field Navier--Stokes vesicle-fluid interaction model,
\textit{J. Math. Fluid Mech.}, \textbf{14} (2012), 177--195.

\bibitem{LowengrubAllardAland2016}
J. Lowengrub, J. Allard, and S. Aland,
Numerical simulation of endocytosis: Viscous flow driven by membranes with non-uniformly distributed curvature-inducing molecules,
\textit{J. Comput. Phys.}, \textbf{309} (2016), 112--128.

\bibitem{LowengrubRatzVoigt2009}
J. Lowengrub, A. R\"{a}tz, and A. Voigt,
Phase-field modeling of the dynamics of multicomponent vesicles: Spinodal decomposition, coarsening, budding and fission,
\textit{Phys. Rev. E}, \textbf{79} (2009), 0311926, 13 pp.

\bibitem{M1}
A. Miranville,
Sixth-order Cahn--Hilliard equations with singular nonlinear terms,
{\it Appl. Anal.}, {\bf 94} (2015), 2133--2146.

\bibitem{M2}
A. Miranville,
On the phase-field-crystal model with logarithmic nonlinear terms,
{\it RACSAM}, {\bf 110} (2016), 145--157.

\bibitem{PZ1}
I. Paw{\l}ow and W. Zajaczkowski,
A sixth order Cahn--Hilliard type equation arising in oil-water-surfactant mixtures,
{\it Comm. Pure Appl. Anal.}, {\bf 10} (2011), 1823--1847.

\bibitem{PZ2}
I. Paw{\l}ow and W. Zajaczkowski,
On a class of sixth order viscous Cahn--Hilliard type equations,
{\it Discrete Contin. Dyn. Syst. Ser. S}, {\bf 6} (2013), 517--546.

\bibitem{PS}
A. Poiatti and A. Signori,
Regularity results and optimal velocity control of the convective nonlocal Cahn--Hilliard equation in 3D,
{\it ESAIM Control Optim. Calc. Var.}, {\bf 30} (2024),
Article Number 21, 36 pages.

\bibitem{PromislowWetton2009}
K. Promislow and B. Wetton,
PEM fuel cells: A mathematical overview,
\textit{SIAM J. Appl. Math.}, \textbf{70} (2009), 369--409.

\bibitem{PQ}
K. Promislow and Q. Wu,
Existence of pearled patterns in the planar functionalized Cahn--Hilliard equation,
{\it J. Differential Equations}, {\bf 259} (2015), 3298--3343.

\bibitem{PromislowWu2017}
K. Promislow and Q. Wu,
Existence, bifurcation, and geometric evolution of quasi-bilayers in the multicomponent functionalized Cahn–Hilliard equation,
\textit{J. Math. Biol.}, \textbf{75} (2017), 443--489.

\bibitem{RoccaSprekels2015}
E. Rocca and J. Sprekels,
Optimal distributed control of a nonlocal convective Cahn--Hilliard equation by the velocity in three dimensions,
\textit{SIAM J. Control Optim.}, \textbf{53} (2015), 1654--1680.

\bibitem{SP}
G. Schimperna and I. Paw{\l}ow,
A Cahn--Hilliard equation with singular diffusion,
{\it J. Differential Equations}, {\bf 254} (2013), 779--803.

\bibitem{SchimpernaWu2020}
G. Schimperna and H. Wu,
On a class of sixth-order Cahn--Hilliard-type equations with logarithmic potential,
\textit{SIAM J. Math. Anal.}, \textbf{52} (2020), 5155--5195.

\bibitem{SeifertLipowsky1995}
U. Seifert and R. Lipowsky,
Morphology of Vesicles,
in \textit{Handbook of Biological Physics}, Vol. 1 (1995), 405--462.

\bibitem{Simon}
J. Simon,
Compact sets in the space $L^p(0,T;B)$.
{\it Ann. Math. Pura. Appl.}, {\bf 146} (1987), 65--96.

\bibitem{Sohr}
H. Sohr,
{\it The Navier--Stokes Equations: An Elementary Functional Analytic Approach},
Birkh\"auser Advanced Texts, Springer, Basel (2010).

\bibitem{SprekelsWu2021}
J. Sprekels and H. Wu,
Optimal distributed control of a Cahn--Hilliard--Darcy system with mass sources,
\textit{Appl. Math. Optim.}, \textbf{83} (2021), 489--530.

\bibitem{TorabiLowengrubVoigtWise2009}
S. Torabi, J. Lowengrub, A. Voigt, and S. Wise,
A new phase-field model for strongly anisotropic systems,
\textit{Proc. R. Soc. Lond. Ser. A Math. Phys. Eng. Sci.}, \textbf{465} (2009), 1337--1359.

\bibitem{Tr10}
F. Tr\"{o}ltzsch,
\textit{Optimal Control of Partial Differential Equations. Theory, Methods and Applications},
Graduate Studies in Mathematics, Vol. 112. AMS, Providence (2010).

\bibitem{Wang2008}
X.-Q. Wang,
Asymptotic analysis of phase-field formulation of bending elasticity models,
\textit{SIAM J. Math. Anal.}, \textbf{39} (2008), 1367--1401.

\bibitem{WW}
C. Wang and S. Wise,
Global smooth solutions of the three-dimensional modified phase field crystal equation,
{\it Methods Appl. Anal.}, {\bf 17} (2010), 191--211.

\bibitem{WuXu2013}
H. Wu and X. Xu,
Strong solutions, global regularity, and stability of a hydrodynamic system modeling vesicle and fluid interactions,
\textit{SIAM J. Math. Anal.}, \textbf{45} (2013), 181--214.

\bibitem{WY}
H. Wu and Y.-C. Yang,
Well-posedness of a hydrodynamic phase-field system for functionalized membrane-fluid interaction.
{\it Discrete Contin. Dyn. Syst. Ser. S}, {\bf 15} (2022), 2345--2389.

\bibitem{YueFengLiuShen2004}
P.-T. Yue, J.-J. Feng, C. Liu, and J. Shen,
A diffuse-interface method for simulating two-phase flows of complex fluids,
\textit{J. Fluid Mech.}, \textbf{515} (2004), 293--317.

\End{thebibliography}

\End{document}
